\documentclass[final,legno]{siamltex}

\usepackage{subfigure}

\usepackage{latexsym}
\usepackage{amsmath}
\usepackage{graphicx}
\usepackage{times}
\usepackage{epic,eepic}
\usepackage{xcolor}
\usepackage{setspace}
\usepackage{amsmath, amssymb}
\usepackage{colortbl}
\usepackage{multirow}
\usepackage{multicol}

\newcommand{\eric}[1]{{\color{black}#1}}

\begin{document}

\title{Mixed Generalized Multiscale
Finite Element Methods and Applications}

\author{
Eric T. Chung \thanks{Department of Mathematics,
The Chinese University of Hong Kong (CUHK), Hong Kong SAR. Email: {\tt tschung@math.cuhk.edu.hk}.
The research of Eric Chung is supported by Hong Kong RGC General Research Fund (Project 400411).}
\and
Yalchin Efendiev \thanks{Department of Mathematics \& Institute for Scientific Computation (ISC),
Texas A\&M University,
College Station, Texas, USA
and Center for Numerical Porous Media (NumPor),
King Abdullah University of Science and Technology (KAUST),
Thuwal 23955-6900, Kingdom of Saudi Arabia. Email: {\tt efendiev@math.tamu.edu}.}
\and
Chak Shing Lee \thanks{Department of Mathematics \& Institute for Scientific Computation (ISC),
Texas A\&M University,
College Station, Texas, USA. Email: {\tt cslee@math.tamu.edu}.}
}

\maketitle

%\author{\textbf{Eric Chung}$^{1}$}
%\author{\textbf{Yalchin Efendiev}$^{3,2}$}

%\author{\textbf{Chak Shing Lee}$^{2}$}

%\address{$^{1}$ Department of Mathematics \\
%The Chinese University of Hong Kong (CUHK) \\
%Hong Kong}

%\address{$^{2}$ Department of Mathematics \& Institute for Scientific Computation (ISC) \\
%Texas A\&M University \\
%College Station, Texas, USA}

%\address{$^{3}$ Center for Numerical Porous Media (NumPor) \\
%King Abdullah University of Science and Technology (KAUST) \\
%Thuwal 23955-6900, Kingdom of Saudi Arabia.}

%\cortext[cor1]{Email address: efendiev@math.tamu.edu}

\begin{abstract}

In this paper, we present a mixed Generalized Multiscale Finite
Element Method (GMsFEM) for solving flow in heterogeneous media. Our
approach constructs multiscale basis functions following a GMsFEM
framework and couples these basis functions using a mixed finite
element method, which allows us to obtain a mass conservative
velocity field.
To construct multiscale basis functions for each coarse edge, we
design a snapshot space that consists of fine-scale velocity fields
supported in a union of two coarse regions that share the common
interface. The snapshot vectors have zero Neumann boundary
conditions on the outer boundaries and we prescribe their values on
the common interface. We describe several spectral decompositions in
the snapshot space motivated by the analysis.
In the paper, we also study oversampling approaches that enhance the accuracy
of mixed GMsFEM. A main idea of oversampling techniques is to introduce
a small dimensional snapshot space.
We present numerical results for two-phase flow and transport,
without updating basis functions in time. Our numerical results show
that one can achieve good accuracy with a few basis functions per
coarse edge if one selects appropriate offline spaces.

\end{abstract}


%\end{frontmatter}

%==============================
\section{Introduction}
%==============================

In many applications, one encounters multiple scales and high
contrast. For example in subsurface applications, the media
properties, such as permeability,
 have multiple
scales and features, such as fractures and shale layers, which have
thickness that are much smaller than the domain size. The solution
techniques for such problems require some type of model reduction
that allows reducing the degrees of freedom. Investigators developed
many approaches, such as upscaling and multiscale methods,
 %\cite{aarnes04,apwy07,eh09,hw97,hughes98,jennylt03}
for this purpose. In these approaches, a coarse grid with a
reduced number of degrees of freedom represents the solution. For
example, in upscaling methods \cite{dur91, weh02}, one typically
upscales the media properties and solves the global problem on a
coarse grid. In multiscale methods
\cite{aarnes04,apwy07,eh09,hw97,jennylt03,ChungLeung,Wave,Reduced,Gibson}, one constructs
multiscale basis functions and solves the problem on a coarse grid.

In this paper, our interest is in developing multiscale mixed
methods for flows in heterogeneous media. Some earlier pioneering
works in multiscale mixed finite element methods are reported in
\cite{ch03, arb00, arbogast02} (see also \cite{apwy07} for mortar
multiscale construction). The motivation for using a mixed finite element
framework is to preserve mass conservation, a property that
is very important for flow problems. In multiscale methods, mixed
methods also provide rapid evaluation of the fine-scale conservative
velocity without doing a post-processing step. The main idea of
mixed multiscale finite element methods is to compute multiscale
basis function for each edge supported in two coarse blocks that
share this common edge (see {\sc Fig}.~\ref{fig:coarse_grid} for
illustration). Computing one basis function per coarse edge using
local or global information limited most previous approaches. As
discussed in our recent paper \cite{egh12, ge09_2,adaptive} (in the framework
of continuous Galerkin approach), one basis function per edge is not
sufficient to capture many disconnected multiscale features; and
therefore, one needs a systematic procedure for enriching the coarse
space. Here, we follow the framework of the Generalized Multiscale
Finite Element Method introduced in \cite{egh12}.

The main idea of the GMsFEM is to divide the computation into
offline and online stages. In the offline stage, we construct (1) a
snapshot space and (2) the offline space via spectral decomposition
of the snapshot space.
 The main concept of constructing the snapshot space
is that the snapshot vectors preserve some essential properties of
the solution and provide a good approximation space. The main idea
of the offline space is that it gives a good solution approximation
with fewer basis functions.
%in the space of snapshot space via local spectral problems.
The main difficulty in the GMsFEM is in finding an appropriate
snapshot space and the local spectral decomposition of the snapshot
space that can give a good approximation of the solution with a
fewer basis functions. We address these issues in this paper.

We construct the velocity field for the snapshot and offline spaces
for the mixed GMsFEM. For approximating the pressure field, we use
piecewise constant basis functions (cf. \cite{ch03, aarnes04,
arbogast02}). The snapshot solutions consist of local solutions with
a unit-variable flux, chosen at all possible locations, on the
internal boundary. To construct the offline space, we consider
several spectral problems based on the analysis and provide a
justification.
%These local spectral problems contain jumps of the
%solution.
We also provide an alternative derivation of the local
eigenvalues problem, where each next vector in the offline space is
furthest from the space of previously selected offline vectors.
%This algorithm can provide a possibly cheaper procedure b
%in future for calculating the offline vectors.
Some of the advantages of the mixed GMsFEM (over the continuous
Galerkin GMsFEM) are the following: (1) no need for partition of
unity; (2) mass conservative and useful for flow and transport.
%(3) it recovers homogenization based mixed MsFEM when there is a scale
%separation.
We also study oversampling techniques (cf. \cite{ehw99})
 by constructing snapshot vectors as the local
solutions in larger regions that contain the interfaces of two adjacent
coarse blocks.
 This allows obtaining a much smaller dimensional
snapshot space and can help to improve the accuracy of the mixed GMsFEM.
Oversampling technique can be particularly helpful for problems
with scale separation. This is because by
taking the restriction of the local solutions in larger domains
in the interior, we avoid the pollution
effects near the boundaries.

We present numerical results for various heterogeneous permeability
fields and show that we can approximate the solution using only a
few basis functions. We use our approach for solving two-phase flow
and transport equations. In two-phase flow and transport, we solve
the flow equation for each time step and employ the fine-scale
velocity to advance the saturation front. In our numerical results,
we solve the flow equation without modifying multiscale basis
functions, i.e., we solve the flow equation on a coarse grid. We
show how adding a few extra basis functions can improve the
prediction accuracy.

We organize the paper as follows. In Section \ref{sec:prelim},
 we present a basic model problem, fine-scale discretization, and the
definitions of coarse and fine grids. In Section
\ref{sec:basisconstruction}, we describe the construction of the snapshot
and offline spaces. We devote Section \ref{sec:analysis} to
analyzing the mixed GMsFEM.
The use of oversampling techniques for the snapshot space is presented in Section \ref{sec:ovs}.
We present numerical results in Section
\ref{sec:numresults}.
The paper ends with a conclusion.

\section{Preliminaries}
\label{sec:prelim}

We consider the following high-contrast flow problem in mixed
formulation
\begin{eqnarray}
\label{pv-system}
\begin{split}
 \kappa^{-1} v+\nabla
p&=0 \quad \text{in} \quad D, \\
\text{div} ( v) &=f \quad \text{in} \quad D,
\end{split}
\end{eqnarray}
with non-homogeneous Neumann boundary condition $v\cdot n = g$ on
$\partial D$, where $\kappa$ is a given high-contrast heterogeneous
permeability field, $D$ is the computational domain, and $n$ is the
outward unit-normal vector on $\partial D$.
%Problem
%(\ref{pv-system}) is required to solve multiple times, and therefore
%some efficient numerical methods are needed.
%There are in literature upscaling methods \cite{dur91}
%and multiscale methods \cite{aarnes04,apwy07,eh09,hw97,jennylt03}
%that provide some efficient solution strategies for the problem (\ref{pv-system%}).
%For media with complex heterogeneities, these methods are not able to produce a%ccurate solutions.
%It is thus the purpose of this paper to develop some enriched multiscale spaces
%which give accurate solutions to (\ref{pv-system}) with few degrees of freedoms.

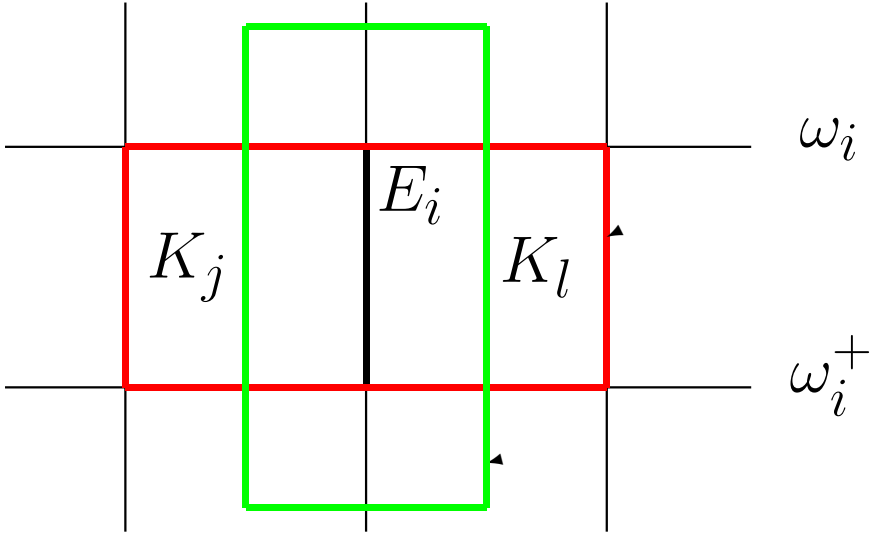
\begin{figure} \setlength{\unitlength}{1.6mm}
\begin{center}
\begin{picture}(80,43)(-5,0)
\linethickness{0.05mm} \color{black!50!white}
\multiput(0,-2)(2,0){32}%
{\dashline[20]{.5}(0,0)(0,44)}
\multiput(0,-2)(0,2){23}%
{\dashline[20]{.5}(0,0)(62,0)} \color{black}\linethickness{0.3mm}
\multiput(10,-2)(20,0){3}%
{\line(0,1){44}}
\multiput(0,10)(0,20){2}%
{\line(1,0){62}} \put(12,18){\colorbox{white}{\makebox(5,4){\Huge
$K_j$}}} \put(41,18){\colorbox{white}{\makebox(5,4){\Huge $K_l$}}}
\put(31,24){\colorbox{white}{\makebox(4,4){\Huge $E_i$}}}
\put(66,9){\colorbox{white}{\makebox(4,4){\Huge $\omega^+_i$}}}
\put(30,10){\linethickness{0.8mm}\line(0,1){20}}
\put(10,10){\linethickness{0.8mm}\color{red}\line(1,0){40}}
\put(10,30){\linethickness{0.8mm}\color{red}\line(1,0){40}}
\put(10,10){\linethickness{0.8mm}\color{red}\line(0,1){20}}
\put(50,10){\linethickness{0.8mm}\color{red}\line(0,1){20}}
\put(65,30){\thicklines\vector(-2,-1){15}} \put(66,30){\Huge
$\omega_i$}
\put(65,10){\thicklines\vector(-4,-1){25}}%\put(35.5,32.2){\thicklines\vector(-2,-4){5.3}}
\linethickness{0.8mm}\color{green}%\put(20,0){\dashline[40]{2}(20,0)(40,0)}
\put(20,0){\line(1,0){20}} \put(20,40){\line(1,0){20}}
\put(20,0){\line(0,1){40}} \put(40,0){\line(0,1){40}}
\end{picture}
\end{center}
\caption{An example of a neighborhood $\omega_i = K_j\cup K_l$ and
$\omega_i^+$ associated with the coarse edge $E_i$.}
\label{fig:coarse_grid}
\end{figure}

In the mixed GMsFEM considered in this paper, we construct the basis
functions for the velocity field, $v=-\kappa\nabla p$. For the
pressure $p$, we will use piecewise constant approximations. To
describe the general solution framework for the model problem
(\ref{pv-system}), we first introduce the notion of fine and coarse
grids. We let $\mathcal{T}^H$ be a usual conforming partition of the
computational domain $D$ into finite elements (triangles,
quadrilaterals, tetrahedrals, etc.), called coarse-grid blocks,
where $H>0$ is the coarse mesh size. We refer to this partition as
the coarse grid and assume that each coarse-grid block is
partitioned into a connected union of fine-grid blocks, which are
conforming across coarse-grid edges. The fine grid partition will be
denoted by $\mathcal{T}^h$, which by definition is a refinement of
$\mathcal{T}^H$. We use $\mathcal{E}^H := \bigcup_{i=1}^{N_e}
\{E_i\}$ (where $N_e$ is the number of coarse edges) to denote the
set of all edges of the coarse mesh $\mathcal{T}^H$, and
$\mathcal{E}^H_0$ to denote the set of all interior coarse edges. We
also define the coarse neighborhood $\omega_i$ corresponding to the
coarse edge $E_i$ as the union of all coarse-grid blocks having the
edge $E_i$, namely,
\begin{equation} \label{neighborhood}
\omega_i=\bigcup\{ K_j\in\mathcal{T}^H; ~~~ E_i\in \partial K_j\}.
\end{equation}
See {\sc Fig.}~\ref{fig:coarse_grid} for an example of a coarse
neighborhood, where the coarse-grid edges are denoted by solid lines
and the fine-grid edges are denoted by dash lines.

%The pressure $p$ is approximated in the space
Let $Q_H$ be the space of piecewise constant functions with respect
to the coarse grid $\mathcal{T}^H$. The approximation of the
pressure $p$ will be obtained in this space. On the other hand, a
set of multiscale basis functions for the velocity field $v$ are
defined for each coarse edge $E_i \in\mathcal{E}^H$ and these basis
functions are supported in the coarse neighborhood $\omega_i$
corresponding to the coarse edge $E_i$. Specifically, to obtain a
basis function for a coarse edge $E_i$, we will solve a local
problem in the coarse neighborhood $\omega_i$ with a given normal
velocity on $E_i$ and zero normal velocity on the boundary
$\partial\omega_i$. Notice that we can use multiple basis functions
for each coarse edge $E_i$ by using various choices of normal
velocity on $E_i$. Let $\{ \Psi_j \}$ be the set of multiscale basis
functions for the edge $E_i$. We define the multiscale space for the
velocity field $v$ as the linear span of all local basis functions
which is denoted as
\begin{equation*}
  {V}_H=\bigoplus_{\mathcal{E}^H}\{ \Psi_i\}.
\end{equation*}
We also define $V_H^0 = V_H \cap \{ v \in V_H \, : \, v \cdot n = 0
\text{ on } \partial D \}$ as a subspace of $V_H$ consisting of
vector fields with zero normal component on $\partial D$; that is,
\begin{equation*}
  {V}^0_H=\bigoplus_{\mathcal{E}_0^H}\{ \Psi_i\}.
\end{equation*}
Given the above spaces, the mixed GMsFEM is to find $( v_H,p_H)\in
{V}_H\times Q_H$ such that
\begin{eqnarray}
\label{approximation-problem}
\begin{split}
\int_D \kappa^{-1} v_H\cdot w_H  - \int_D \text{div} (w_H)p_H
&=0,\quad
\forall w_H\in  {V}_H^0, \\
\int_D \text{div} (v_H) q_H  &=\int_D f q_H, \quad \forall q_H \in
Q_H,
\end{split}
\end{eqnarray}
where $v_H \cdot n = g_H$ on $\partial D$, and for each coarse edge
$E_i\in\partial D$, we have
\begin{equation*}
\int_{E_i} (g_H - g) \Psi_j\cdot n = 0
\end{equation*}
for all basis functions $\Psi_j$ corresponding to the edge $E_i$.

%$v_h\cdot {n}=g_h $ on $\partial D$,
%where $g_h={  g}_{0,h} {  n}$ on $\partial D$
%and ${  g}_{0,h}=\sum_{e\in\{ \partial K\bigcap \partial\Omega, K\in {\cal T}_h \}}
%(\int_e g ds){  N}_e$, ${  N}_e\in {  V}_h$, is corresponding basis
%function to edge $e$,

In addition, we let $V_h \times Q_h$ be the standard lowest-order
Raviart-Thomas space for the approximation of (\ref{pv-system}) on
the fine grid $\mathcal{T}^h$. Then, the fine-grid solution
$(v_h,p_h)$ satisfies
\begin{eqnarray}
\label{fine-problem}
\begin{split}
\int_D \kappa^{-1} v_h\cdot w_h  - \int_D \text{div} (w_h)p_h
&=0,\quad
\forall w_h \in  {V}_h^0, \\
\int_D \text{div} (v_h) q_h  &=\int_D f q_h, \quad \forall q_h \in
Q_h,
\end{split}
\end{eqnarray}
where $v_h \cdot n = g_h$ on $\partial D$ and $V_h^0 = V_h \cap \{ v
\in V_h \, : \, v \cdot n = 0 \text{ on } \partial D \}$. In terms
of matrix representations, the above problem can be written as
\begin{equation}
\begin{split}
M_{\text{fine}} \vec{V}_h + B_{\text{fine}}^T \vec{P}_h &= 0 \\
B_{\text{fine}}\vec{V}_h &= \vec{F}_h,
\end{split}
\label{eq:finematrix}
\end{equation}
where $\vec{V}_h$ and $\vec{P}_h$ are vectors of coefficients in the
expansions of the solutions $v_h$ and $p_h$ in the spaces $V_h$ and
$Q_h$, respectively. We remark that the fine-grid solution
$(v_h,p_h)$ is considered as a reference solution, and we will
compare the accuracy of the multiscale solution $( v_H,p_H)$ against
the fine grid solution. Furthermore, it is easy to see that $Q_H
\subset Q_h$. We will construct the multiscale space $V_H$ so that
$V_H \subset V_h$. Therefore, the mixed GMsFEM can be considered as
a conforming method to approximate the fine-grid solution. In the
next section, we will give the construction of the multiscale basis
functions and the space $V_H$.

\section{The construction of multiscale basis functions}
\label{sec:basisconstruction}

In this section, we will discuss the construction of the multiscale
space $V_H$ for the approximation of the velocity field. We will
first introduce the snapshot space, which contains an extensive set
of basis functions formed by solutions of local problems with all
possible boundary conditions up to the fine-grid resolution. Then,
we will present a space reduction technique which provides a
systematic way to select the dominant modes in the snapshot space.
This technique is based on a carefully designed local spectral
problem giving a rapidly decaying residual. The resulting reduced
space is obtained by the linear span of these dominant modes and is
called the offline space. Notice that we use the terminology
introduced in \cite{egh12}, where the notion of online space is also
introduced. We emphasize that, since we consider problems without
parameter dependence, the offline space is the same as the online
space.

\subsection{Snapshot space}\label{sec:snap}

In this section, we will define the snapshot space. Essentially, it
is a space containing an extensive set of basis functions which are
solutions of local problems with all possible boundary conditions up
to the fine-grid resolution. Specifically, the functions in the
snapshot space are $\kappa$-harmonic functions of unit-flux
functions. In the following, we explain the detailed constructions.
Let $E_i \in\mathcal{E}^H$ be a coarse edge. We will find
$(v_j^{(i)}, p_j^{(i)})$ by solving the following problem on the
coarse neighborhood $\omega_i$ corresponding to the edge $E_i$
\begin{eqnarray}
\label{MMs-basis-eq}
\begin{split}
 \kappa^{-1} v_j^{(i)}+\nabla
p_j^{(i)} &=0 \quad\quad \text{in} \quad \omega_i, \\
\text{div} ( v_j^{(i)}) &= \alpha_j^{(i)} \quad \text{in} \quad
\omega_i,
\end{split}
\end{eqnarray}
subject to the boundary condition $v_j^{(i)}\cdot n_i = 0$ on
$\partial\omega_i$, where $n_i$ denotes the outward unit-normal
vector on $\partial\omega_i$. One key feature of our proposed
approach is that the above problem (\ref{MMs-basis-eq}) will be
solved separately in the coarse-grid blocks forming $\omega_i$ (see
Figure \ref{fig:coarse_grid} for illustration). Therefore, we will
need an extra boundary condition on $E_i$, which is discussed below.
Notice that the coarse edge $E_i$ can be written as a union of
fine-grid edges, namely, $E_i = \bigcup_{j=1}^{J_i} e_j$, where
$J_i$ is the total number of fine-grid edges on $E_i$ and $e_j$
denotes a fine-grid edge. Let $\delta_j^{(i)}$ be a piecewise
constant function defined on $E_i$ with respect to the fine-grid
such that it has value $1$ on $e_j$ and value $0$ on the other
fine-grid edges; that is,
\begin{equation*}
\delta_j^{(i)} = \left\{
\begin{array}{ll}
1, &\mbox{on $e_j$,}\\
0, &\mbox{on other fine grid edges on $E_i$,}
\end{array}
\right. \quad\quad j = 1,2,\cdots, J_i.
\end{equation*}
The remaining boundary condition on the coarse edge $E_i$ for the
local problem (\ref{MMs-basis-eq}) is then taken as
\begin{equation}
v_j^{(i)} \cdot m_i = \delta_j^{(i)} \quad \mbox{on $E_j$,}
\label{eq:BC}
\end{equation}
%\begin{eqnarray}
%\label{MMs-basis-eq}
%\begin{split}
%\text{div} (\kappa(x)\nabla w_j^{(i)}) &=& C \quad \text{in} \quad  \omega_i\\
 %\kappa(x)\nabla w_j^{(i)}\cdot {n} &=& \left\{
%\begin{array}{ll}
%\delta_j &\mbox{on $e_j$}\\
%0 &\mbox{else,}
%\end{array}
%\right.
%\end{split}
%\end{eqnarray}
where $m_i$ is a fixed unit-normal vector on $E_i$.
%$\delta_i$ is a fine-grid piecewise constant function that is
%$1$ on a fine-grid edge $e_i$ (see Figure for illustration) and
%$C(x)$ is a piecewise
 %constant function that is constant
%in each coarse block $K$ with the value chosen such that the mass
%balance is satisfied.
We remark that the constant $\alpha_j^{(i)}$ in (\ref{MMs-basis-eq})
is chosen so that the compatibility condition $\int_{K_l}
\alpha_j^{(i)} = \int_{E_i} \delta_j^{(i)}$ is satisfied, for all
$K_l \subset \omega_i$.
%We also remark that the above problem (\ref{MMs-basis-eq})
%is solved separately on the fine grid partition of $K_l$ and $K_j$,
%which are the two coarse grid blocks forming $\omega_i$.
We also remark that, since $v_j^{(i)} \cdot n_i = 0$ on the boundary
of $\omega_i$, the vector field $v_j^{(i)}$ can be extended to the
rest of the domain $D$ by defining $v_j^{(i)}=0$ outside $\omega_i$.
Furthermore, the above local problem (\ref{MMs-basis-eq}) can be
solved numerically on the underlying fine grid of $\omega_i$ by the
lowest-order Raviart-Thomas element, so that the resulting velocity
$v_j^{(i)} \in V_h$ (for simplicity, we keep the same notation for
the discrete solution $v_j^{(i)}$).

The collection of the solutions of the above local problems
generates the snapshot space. We let $\Psi^{i,\text{snap}}_j :=
v_j^{(i)}$ be the snapshot fields and define the snapshot space
$V_{\text{snap}}$ by
\begin{equation*}
V_{\text{snap}} = \text{span}\{ \Psi_{j}^{ i,\text{snap}}:~~~1\leq j
\leq J_i, \quad 1\leq i\leq N_e \}.
\end{equation*}
%{\bf{Throughout, we will not distinguish between
%the fine-grid solutions of local problems and their continuous counterparts.}}
To simplify notation, we will use the following single-index
notation
\begin{equation*}
V_{\text{snap}} = \text{span}\{ \Psi_{i}^{ \text{snap}}:~~~ 1\leq
i\leq M_{\text{snap}} \},
\end{equation*}
where $M_{\text{snap}} = \sum_{i=1}^{N_e} J_i$ is the total number
of snapshot fields.

%Next, we define $\Psi_i^{\text{snap}}=-\kappa\nabla w_i$ as snapshot vectors.
%Via a spectral decomposition of the snapshot space, we will identify
%the offline space $\Psi_i$'s. This will be presented in the next section.

%We denote the snapshot space in each $\omega$ by
%$V_{\text{snap}}$,
%$$
%V_{\text{snap}} = \text{span}\{ \Psi_{j}^{ \text{snap}}:~~~1\leq j \leq J \},
%$$
%where  $J$ denotes the number of snapshots.

Notice that each $\Psi^{\text{snap}}_i$ is represented on the fine
grid by the basis functions in $V_h$. Therefore, each
$\Psi^{\text{snap}}_i$ can be represented by a vector
$\psi_i^{\text{snap}}$ containing the coefficients in the expansion
of $\Psi^{\text{snap}}_i$ in the fine-grid basis functions. Then, we
define
$$
R_{\text{snap}} = \left[ \psi_{1}^{\text{snap}}, \ldots,
\psi_{M_{\text{snap}}}^{\text{snap}} \right],
$$
which maps from the coarse space to the fine space.

\subsection{Offline space}\label{sec:offline}

Following the general framework of \cite{egh12}, we will perform a
space reduction on the snapshot space through the use of some local
spectral problems. The reduced space is called the offline space.
The purpose of this is to determine the important modes in the
snapshot space and to obtain a smaller space for approximating the
solution. In the general setting, we consider the spectral problem
of finding a real number $\lambda$ and a vector field $v\in
V_{\text{snap}}$ such that
\begin{equation}
\label{eq:generalspec} a(v,w) = \lambda s(v,w), \quad\quad \forall
w\in V_{\text{snap}},
\end{equation}
where $a(v,w)$ and $s(v,w)$ are symmetric positive definite bilinear
forms defined on $V_{\text{snap}}\times V_{\text{snap}}$. We
consider $s(v,w)$ as an inner product on $V_{\text{snap}}$ and
define a linear operator $\mathcal{A}: V_{\text{snap}}\rightarrow
V_{\text{snap}}$ by
\begin{equation*}
s(\mathcal{A}v,w) = a(v,w).
\end{equation*}
We assume that the operator $\mathcal{A}$ has rapidly decaying
eigenvalues.
%is compact so that
%the eigenvalues of the above spectral problem (\ref{eq:generalspec})
%are rapidly increasing.
Note that one can take $\mathcal{A}$ to be a compact operator.

In practice, solving the above global spectral problem
(\ref{eq:generalspec}) is inefficient. Therefore, the dimension
reduction and the construction of the offline space are performed
locally. In particular, the above spectral problem is solved for
each coarse neighborhood $\omega_i$ corresponding to the coarse edge
$E_i$. We let $V_{\text{snap}}^{(i)}$ be the snapshot space
corresponding to the coarse edge $E_i$, which is defined by
\begin{equation*}
V^{(i)}_{\text{snap}} = \text{span}\{ \Psi_{j}^{
i,\text{snap}}:~~~1\leq j \leq J_i \}.
\end{equation*}
The local spectral problem is: find a real number $\lambda \geq 0$
and a function $v\in V_{\text{snap}}^{(i)}$ such that
\begin{equation}
\label{eq:localspec} a_i(v,w) = \lambda s_i(v,w), \quad\quad \forall
w\in V_{\text{snap}}^{(i)}.
\end{equation}
We will consider two different choices of local spectral problems.
%For both of these choices, we take
%\begin{equation*}
%s_i(v,w) = \int_{\omega_i} \kappa^{-1} v\cdot w.
%\end{equation*}
%Regarding the bilinear form $a_i(v,w)$, we consider the following
%two options, which are shown numerically to be promising ways to
%obtain reduced spaces.
One can possibly use oversampling ideas
\cite{hw97, ehw99, bl11, eglp13}
 to achieve a better convergence rate (see Section \ref{sec:ovs}).

\noindent {\bf Spectral problem 1}: We take
\begin{equation}
a_i(v,w) = \int_{E_i} \kappa^{-1} (v\cdot m_i) (w\cdot m_i),
\quad
s_i(v,w) = \int_{\omega_i} \kappa^{-1} v\cdot w + \int_{\omega_i} \text{div}(v) \, \text{div}(w),
\label{eq:pod1}
\end{equation}
where we recall that $m_i$ is a fixed unit-normal on the coarse edge
$E_i$.

\noindent {\bf Spectral problem 2}: We take
\begin{equation}
a_i(v,w) = \int_{\omega_i} \kappa^{-1} v\cdot w, \quad s_i(v,w) = \int_{E_i} [p_v] [p_w],
\label{eq:pod2}
\end{equation}
where $(v,p_v)$ and $(w,p_w)$ are solutions of the local problem
(\ref{MMs-basis-eq}), and $[p]$ denotes the jump of the function
$p$.

%In fact, the above two spectral problems give the same
%eigenfunctions and the corresponding eigenvalues are reciprocal to each other.
%We illustrate this for an interior edge $E_i$.
%Writing $\omega_i = K_j \cup K_l$, we have
%\begin{equation*}
%s_i(v,w) = \int_{K_j} \kappa^{-1} v\cdot w + \int_{K_l} \kappa^{-1} v\cdot w.
%\end{equation*}

In the following, we will focus our discussions on spectral problem 1.
For spectral problem 2, we will only report
its performance in Section \ref{sec:numresults}
to show that it is also a promising way to obtain a reduced space.

Assume that the eigenvalues of (\ref{eq:localspec}) are arranged in
increasing order
\begin{equation}
\label{eq:eigenorder} \lambda_1^{(i)} < \lambda_2^{(i)} < \cdots  <
\lambda_{J_i}^{(i)},
\end{equation}
where $\lambda_k^{(i)}$ denotes the $k$-th eigenvalue for the coarse
neighborhood $\omega_i$. The corresponding eigenvectors are denoted
by $Z_k^{(i)} = (Z_{kj}^{(i)})_{j=1}^{J_i}$, where $Z_{kj}^{(i)}$ is
the $j$-th component of the vector $Z_k^{(i)}$. We will use the
first $l_i$ eigenfunctions to form the offline space. We remark that
we assume the eigenvalues are strictly increasing (here, we refer to
the inverse of $\mathcal{A}$, cf. (\ref{eq:generalspec})) only to
simplify the discussion. In practice, if there are multiple
eigenvectors corresponding to a specific eigenvalue, then we will
take all these eigenvectors to be part of the basis functions when
the corresponding eigenvalue is selected. Using the eigenfunctions,
offline basis functions can be constructed as
\begin{equation*}
\Psi_k^{i,\text{off}} = \sum_{j=1}^{J_i} Z_{kj}^{(i)}
\Psi_j^{i,\text{snap}}, \quad\quad k=1,2,\cdots, l_i.
\end{equation*}
The global offline space is then
\begin{equation*}
V_{\text{off}} = \text{span}\{ \Psi_{k}^{ i,\text{off}}:~~~1\leq k
\leq l_i, \quad 1\leq i\leq N_e \}.
\end{equation*}
To simplify notation, we will use the following single-index
notation
\begin{equation*}
V_{\text{off}} = \text{span}\{ \Psi_{k}^{ \text{off}}:~~~ 1\leq
k\leq M_{\text{off}} \},
\end{equation*}
where $M_{\text{off}} = \sum_{i=1}^{N_e} l_i$ is the total number of
offline basis functions. This space will be used as the
approximation space for the velocity; that is, $V_H =
V_{\text{off}}$ in the GMsFEM system (\ref{approximation-problem}).
Furthermore, we define $V_{\text{off}}^0$ as the restriction of
$V_{\text{off}}$ formed by the linear span of all basis functions
$\Psi_k^{\text{off}}$ corresponding to interior coarse edges only.
Thus, all vectors in $V_{\text{off}}^0$ have zero normal component
on the global domain boundary $\partial D$.

%In order to construct the offline space $V_{\text{off}}^\omega$, we perform a dimension reduction of the space of snapshots using an auxiliary spectral decomposition. We seek a subspace of the snapshot space such that it can approximate any element of the snapshot space in the appropriate sense defined via auxiliary bilinear forms.
%
In term of matrix representations, the above eigenvalue problem
(\ref{eq:localspec}) can be expressed as
\begin{equation} \label{offeig}
A_{\text{snap}}^{(i)} Z_k^{(i)} = \lambda_k^{(i)}
S_{\text{snap}}^{(i)} Z_k^{(i)},
\end{equation}
where
\begin{equation*}
 \displaystyle A_{\text{snap}}^{(i)} = [(A_{\text{snap}}^{(i)})_{mn}]
 = a_i( \Psi_m^{i,\text{snap}} ,  \Psi_n^{i,\text{snap}}) = R_{\text{snap}}^T A_{\text{fine}}^{(i)} R_{\text{snap}}
 \end{equation*}
 \begin{center}
 and
 \end{center}
 \begin{equation*}
 \displaystyle S_{\text{snap}}^{(i)}
 = [(S_{\text{snap}}^{(i)})_{mn}] = s_i( \Psi_m^{i,\text{snap}} ,  \Psi_n^{i,\text{snap}}) = R_{\text{snap}}^T S_{\text{fine}}^{(i)} R_{\text{snap}}.
\end{equation*}
 We note that $A_{\text{fine}}^{(i)}$ and $S_{\text{fine}}^{(i)}$ denote analogous fine-scale matrices that use fine-grid basis functions.
 %To generate the offline space we then choose the smallest $M_{\text{off}}$ eigenvalues  and form the corresponding eigenvectors in the space of snapshots by setting
%$\psi_k^{\text{off}} = \sum_j \Psi_{kj}^{\text{off}} \psi_j^{\text{snap}}$ (for $k=1,\ldots, M_{\text{off}}$), where $\Psi_{kj}^{\text{off}}$ are the coordinates of the vector $\Psi_{k}^{\text{off}}$. We then create the offline matrix $$
%R_{\text{off}} = \left[ \psi_{1}^{\text{off}}, \ldots, \psi_{M_{\text{off}}}^{\text{off}} \right]
%$$
%to be used in the online space construction.
Notice that each $\Psi^{\text{off}}_k$ is represented on the fine
grid. Therefore, each $\Psi^{\text{off}}_k$ can be represented by a
vector $\psi_k^{\text{off}}$ containing the coefficients in the
expansion of $\Psi^{\text{off}}_k$ in the fine-grid basis functions.
Then, we define
$$
R_{\text{off}} = \left[ \psi_{1}^{\text{off}}, \ldots,
\psi_{M_{\text{off}}}^{\text{off}} \right],
$$
which maps from the offline space to the fine space. Similar to
(\ref{eq:finematrix}), the GMsFEM system
(\ref{approximation-problem}) can be represented in matrix form as
follows.
\begin{equation}
\begin{split}
R_{\text{off}}^T M_{\text{fine}} R_{\text{off}} \vec{V}_H + R_{\text{off}}^T B_{\text{fine}}^T G_H \vec{P}_H &= 0 \\
G_H^T B_{\text{fine}} R_{\text{off}} \vec{V}_H &= G_H^T \vec{F}_h,
\end{split}
\label{eq:coarsematrix}
\end{equation}
where $G_H$ is the restriction operator from $Q_H$ into $Q_h$, and
$\vec{V}_H$ and $\vec{P}_H$ are vectors of coefficients in the
expansions of the solutions $v_H$ and $p_H$ in the spaces $V_H$ and
$Q_H$, respectively. From (\ref{eq:coarsematrix}), it is easy to see
that implementing the mixed GMsFEM requires the construction of the
fine-grid matrices $M_{\text{fine}}$ and $B_{\text{fine}}$ as well
as the offline matrix $R_{\text{off}}$.

Next, we discuss the eigenvalue behavior which is important for the
method. First, we note that the eigenvalues increase to infinity as
we refine the fine grid. More precisely, the large eigenvalues scale
as the inverse of the fine-scale mesh size. Further,  we note that
the first eigenvector can be chosen to be the multiscale basis function
defined in the mixed MsFEM presented by Chen and Hou \cite{ch03}.
%The second eigenvalue scales
%as $H^{-1/2}$ if the coefficients are bounded and the solution is
%in $W^{1,\infty}$. This indicates that if the problem has a scale
%separation, one will obtain the resonance error in \cite{ch03}
%that will come from the first basis function and
%an additional error that will scale as the coarse mesh size.
%This shows that our approach recovers the classical mixed approach in
%the case of scale separation.

\subsection{Optimization viewpoint of the basis functions}

In this section, we present an optimization viewpoint for the basis
functions obtained by the local spectral problem
(\ref{eq:localspec}). Recall that, for each coarse neighborhood
$\omega_i$, we will solve the spectral problem (\ref{eq:localspec})
to get a sequence of eigenpairs $(\lambda_k^{(i)}, Z_k^{(i)})$. We
will show, by means of an optimization approach, that the
eigenfunction $Z_k^{(i)}$ is furthest away from the space spanned by
the previous eigenvectors $Z_1^{(i)}, \cdots, Z_{k-1}^{(i)}$. Thus,
whenever a new basis function is added, this basis function will
represent an important component in the solution space.

Assume that $k-1$ basis functions, $\phi_1,\cdots, \phi_{k-1}$, are
selected for a specific coarse neighborhood $\omega_i$. Let $W$ be
the space spanned by these functions. To find an additional basis
function, we will find a function $\phi_k$ orthogonal to the space
$W$ and furthest away from the space $W$. To be more specific, we
let $W^{\perp}$ be the orthogonal complement of $W$ with respect to
the inner product defined by the bilinear form $s_i(v,w)$; namely,
\begin{equation*}
W^{\perp} = \Big\{ v \in V_{\text{snap}}^{(i)} \; | \; s_i(v,w)=0,
\; \forall w\in W \Big\}.
\end{equation*}
Then, the function $\phi_k$ is obtained by the following constrained
optimization problem
\begin{equation*}
\begin{split}
& \phi_k = \arg\max_{\phi \in W^{\perp}} s_i(\phi - w,\phi-w), \\
&\text{subject to } \; a_i(\phi,\phi) = 1,
\end{split}
\end{equation*}
for all $w\in W$. By orthogonality, the above problem can be
formulated as
\begin{equation*}
\begin{split}
& \phi_k = \arg\max_{\phi \in W^{\perp}} s_i(\phi,\phi), \\
&\text{subject to } \; s_i(\phi,\phi) = 1.
\end{split}
\end{equation*}
It is well-known that the Euler-Lagrange equation for the above
optimization problem is
\begin{equation*}
\begin{split}
& s_i(\phi_k,\phi) - \mu a_i(\phi_k,\phi) = 0,\quad\forall \phi \in W^{\perp}, \\
& a_i(\phi_k,\phi_k) = 1,
\end{split}
\end{equation*}
where $\mu$ is the Lagrange multiplier. The above condition explains
why we select the eigenfunctions of the spectral problem
(\ref{eq:localspec}) as basis functions.

\subsection{Postprocessing}\label{sec:post}

In this section, we present a postprocessing technique to enhance
the conservation property of the mixed GMsFEM solution. First,
notice that the mixed GMsFEM is conservative on the coarse-grid
level. Specifically, the solution of (\ref{approximation-problem})
satisfies
\begin{equation}
\int_{\partial K} v_H\cdot n  =\int_{K} f
\end{equation}
for every coarse-grid block $K$. This is a direct consequence of the
second equation of (\ref{approximation-problem}) and the fact that
$Q_H$ contains functions that are constant in each coarse block.
%However,
%this coarse level conservation is insufficient for many flow and transport simulations, for which a
%conservative velocity field in the fine grid level is desired.
When $f$ has fine-scale oscillation in some coarse blocks, the
velocity field needs to be postprocessed in these coarse blocks. In
porous media applications, there are only a few coarse blocks where
the sources and sinks are.
%This means that we want to have a discrete
%function $\tilde{v}_H$ defined on the fine grid satisfies the
%condition
In the following, we will construct a postprocessed velocity
$v_h^{\star}$ such that conservation on the fine grid is obtained,
that is,
\begin{equation}
\int_{\partial \tau} (v_h^{\star}\cdot n)  =\int_{\tau} f,
 \;\;\;\; \forall \tau\in\mathcal{T}^h.\label{eq:fine_conserv}
\end{equation}
%In general this does not hold for the solution of our mixed
%generalized multiscale finite element method unless the source term
%$f\in Q_h$. Thus, we need to do a postprocess to the mixed GMsFEM
%velocity solution $v_H$ so that the postprocessed solution
%$\tilde{v}_H$ satisfies (\ref{eq:fine_conserv}).
%We emphasize that the postprocessing is performed locally.
%is simple and local. First, we take the normal trace of
%$\tilde{v}_H$ to be the same as that of $v_H$ on $\partial K_i$.
In particular, for each coarse-grid block $K$, we find
$(v_h^{\star},p_h^{\star}) \in V_h(K)\times Q_h(K)$ such that
$v_h^{\star}\cdot n = v_H\cdot n$ and
%Then we solve the each of the local problems
\begin{eqnarray}
\label{local-problem}
\begin{split}
\int_{K} \kappa^{-1} v_h^{\star} \cdot w_h  - \int_{K}
p_h^{\star}\;\text{div}(w_h)  &=0,\quad
\forall w_h\in  {V}^0_h(K)\\
\int_{K} \text{div}(v_h^{\star}) q_h  &=\int_{K} f q_h, \quad
\forall q_h \in Q_h(K).
\end{split}
\end{eqnarray}
%where ${V}_h^i$ and $Q_h^i$ are restriction of fine scale finite
%element spaces on the coarse grid block $K_i$. Then
%(\ref{eq:fine_conserv}) is true because of the second equation of
%(\ref{local-problem}). Moreover, $\tilde{v}_H$ satisfies the
%compatibility condition
%\begin{equation}
%\int_{\partial K_i} (\tilde{v}_H\cdot n) \;d\sigma =\int_{\partial
%K_i} (v_H\cdot n) \;d\sigma =\int_{K_i} f \;dx.
%\end{equation}
In the single-phase and two-phase flow and transport simulation
experiments below, we will apply this postprocessing technique to
obtain conservative velocity fields on the fine-grid level. We
remark that this postprocessing is only needed in the coarse blocks
where the source term $f$ is non-constant. Therefore, computing the
postprocessed velocity $v_h^{\star}$ is very efficient.
%Therefore, in the flow simulations
%below, we only need to do the postprocess in the coarse blocks which
%are located at the top-left and bottom-right corners of the
%computational domain.

\section{Convergence of the mixed GMsFEM}
\label{sec:analysis}

In this section, we will prove the convergence of the mixed GMsFEM
(\ref{approximation-problem}). The analysis consists of two main
steps. In the first step, we will construct a projection of the
fine-grid velocity field $v_h$ to the snapshot space, and derive an
error estimate for such projection. In the second step, we will
derive an estimate for the difference between the projection of the
fine-grid velocity and the GMsFEM solution. Combining the above two
steps, we obtain an estimate for the difference between the
fine-grid and the GMsFEM solution.

Recall that $(v_h,p_h) \in V_h \times Q_h$ is the fine-grid solution
obtained in (\ref{fine-problem}). We will define a projection
$\widehat{v} \in V_{\text{snap}}$ as follows. Let $K$ be a
coarse-grid block and let $\overline{f} = \frac{1}{|K|} \int_K f$ be
the average value of $f$ over $K$. Then, the restriction of
$\widehat{v}$ on $K$ is obtained by solving the following problem
\begin{eqnarray}
\label{proj-v}
\begin{split}
 \kappa^{-1} \widehat{v} + \nabla
\widehat{p} &=0 \quad\quad \text{in} \quad K, \\
\text{div} ( \widehat{v} ) &= \overline{f} \quad\quad \text{in}
\quad K,
\end{split}
\end{eqnarray}
subject to the following conditions
\begin{equation}
\widehat{v}\cdot n = v_h \cdot n, \; \text{ on}\; \partial K
\quad\text{and}\quad \int_K \widehat{p} = \int_K p_h.
\label{proj-v1}
\end{equation}
We remark that the above problem (\ref{proj-v})-(\ref{proj-v1}) is
solved on the fine grid, and therefore we have $\widehat{v}\in V_h$.
By the construction, we also have $\widehat{v}\in V_{\text{snap}}$.

Now, we introduce some notations for the following analysis. Let
$\Omega$ be an open set. For a scalar function $q\in L^2(\Omega)$,
the $L^2$ norm is $ \| q\|_{L^2(\Omega)}^2 = \int_{\Omega} q^2 $;
and for a vector field $v$, we define the weighted $L^2$ norm
$\|v\|^2_{\kappa^{-1},\Omega} = \int_{\Omega} \kappa^{-1} |v|^2$.
Moreover, the notation $H(\text{div};\Omega;\kappa^{-1})$ denotes
the standard Sobolev space containing vector fields $v$ with $v\in
L^2(\Omega)^2$ and $\text{div}(v)\in L^2(\Omega)$, equipped with
norm $\|v\|_{H(\text{div};\Omega);\kappa^{-1}}^2 = \|
v\|_{\kappa^{-1},\Omega} + \| \text{div}(v) \|_{L^2(\Omega)}^2$. If
$\kappa=1$, we write $H(\text{div};\Omega) =
H(\text{div};\Omega;\kappa^{-1})$. Furthermore, $\alpha \preceq
\beta$ means that there is a uniform constant $C>0$ such that the
two quantities $\alpha$ and $\beta$ satisfy $\alpha \leq C \beta$.

Next, we prove the following estimate for $\widehat{v}$.
\begin{lemma}
\label{lemma1} Let $(v_h,p_h) \in V_h \times Q_h$ be the fine-grid
solution obtained in (\ref{fine-problem}) and $\widehat{v}\in V_h
\cap V_{\text{snap}}$ be the solution of
(\ref{proj-v})-(\ref{proj-v1}). We have
\begin{equation}
\int_D \kappa^{-1} | v_h - \widehat{v}|^2 \preceq
\max_{K\in\mathcal{T}^H} \Big( \kappa^{-1}_{\text{min},K} \Big)
\sum_{i=1}^{N_e} \| f - \overline{f}\|_{L^2(K_i)}^2, \label{bound1}
\end{equation}
where $\kappa_{\text{min},K}$ is the minimum of $\kappa$ over $K$
\end{lemma}

{\it Proof}. Let $K\in\mathcal{T}^H$ be a given coarse-grid block.
First, substracting (\ref{fine-problem}) by the variational form of
(\ref{proj-v}), we have
\begin{eqnarray}
\label{error-eqn}
\begin{split}
\int_K \kappa^{-1} (v_h - \widehat{v}) \cdot w_h  - \int_K
\text{div} (w_h) (p_h - \widehat{p}) =0,\quad
&\forall w_h \in  {V}_h^0(K), \\
\int_K \text{div} (v_h-\widehat{v}) q_h  =\int_K (f-\overline{f})
q_h, \quad &\forall q_h \in Q_h(K),
\end{split}
\end{eqnarray}
where $Q_h(K)$ is the restriction of $Q_h$ on $K$ and $V_h^0(K)$ is
the restriction of $V_h$ on $K$ containing vector fields with zero
normal component on $\partial K$. Taking $w_h = v_h-\widehat{v}$ and
$q_h = p_h - \widehat{p}$ in (\ref{error-eqn}), and summing up the
resulting equations, we have
\begin{equation}
\int_K \kappa^{-1} (v_h - \widehat{v}) \cdot (v_h - \widehat{v}) =
\int_K (f-\overline{f}) (p_h - \widehat{p}). \label{error-eqn-1}
\end{equation}
Recall that the Raviart-Thomas element satisfies the following
inf-sup condition \cite{Brezzi_Fortin_book}:
\begin{equation}
\| q_h \|_{L^2(K)} \preceq \sup_{w_h\in V_h(K)} \frac{ \int_K
\text{div}(w_h) q_h }{ \|w_h\|_{H(\text{div};K)} }, \quad\forall q_h
\in Q_h(K), \label{eq:infsup}
\end{equation}
where $V_h(K)$ is the restriction of $V_h$ on $K$. Using the inf-sup
condition (\ref{eq:infsup}) and the error equation
(\ref{error-eqn}), we have
\begin{equation*}
\| p_h - \widehat{p} \|_{L^2(K)} \preceq
\kappa_{\text{min},K}^{-\frac{1}{2}} \| v_h -
\widehat{v}\|_{\kappa^{-1},K}.
\end{equation*}
Finally, by (\ref{error-eqn-1}), we obtain
\begin{equation*}
 \| v_h - \widehat{v}\|_{\kappa^{-1},K} \preceq \kappa_{\text{min},K}^{-\frac{1}{2}} \| f - \overline{f} \|_{L^2(K)}.
\end{equation*}
Collecting results for all coarse-grid blocks, we obtain the desired
estimate (\ref{bound1}).

\begin{flushright}
$\square$
\end{flushright}

To simplify the notations, we will consider the case with
homogeneous Neumann boundary condition in (\ref{pv-system}). In this
case, the multiscale basis functions are obtained only for interior
coarse edges. We emphasize that the same analysis can be applied to
the non-homogeneous case. Let $N_0$ be the number of interior coarse
edges. For each interior coarse edge $E_i$, we assume that there
exists a basis function $\Psi^{i,\text{off}}_{r_i} \in
V_{\text{off}}^0$, $1\leq r\leq l_i$, such that $\int_{E_i}
\Psi^{i,\text{off}}_{r_i} \cdot m_i \ne 0$. We remark that this is a
reasonable assumption otherwise all basis functions are divergence
free.
%We use the notation $
%We also assume that the index $r_i$ is the largest index with
%the above property.
%We let $\mu^{(i)}$ be the maximum of all eigenvalues $\lambda^{(i)}_r$ with the corresponding basis function satisfies $\int_{E_i} \Psi^{i,\text{off}}_r \cdot m_i \ne 0$.
As a key step in the proof of the main result in Theorem \ref{thm1},
we first prove the following inf-sup condition.

\begin{theorem}
For all $p\in Q_H$, we have
\begin{equation}
\| p \|_{L^2(D)} \preceq  C_{\text{\rm infsup}} \sup_{w\in V^0_{\text{\rm off}}} \frac{
\int_D \text{\rm div}(w) p }{ \| w\|_{H(\text{\rm
div};D;\kappa^{-1})} }, \label{eq:infsup1}
\end{equation}
where $C_{\text{\rm infsup}} = \Big( \max_{1\leq i\leq N_0} \min_{r}
\int_{\omega_i} \kappa^{-1}  \Psi^{i,\text{\rm off}}_{r} \cdot
\Psi^{i,\text{\rm off}}_{r} + 1 \Big)^{\frac{1}{2}}$
and the minimum is taken over all indices $r$ with the property
$\int_{E_i} \Psi^{i,\text{\rm off}}_r \cdot m_i \ne 0$.
\end{theorem}

{\it Proof}. Let $p\in Q_H$. We consider the following Neumann
problem
\begin{equation*}
\begin{split}
\Delta \zeta &= p, \quad \text{ in } D, \\
\frac{\partial \zeta}{\partial n} &= 0, \quad \text{ on } \partial
D.
\end{split}
\end{equation*}
We assume that the solution $\zeta\in H^2(D)$ and we let $\eta =
\nabla \zeta$. Then we will define $w\in V_{\text{off}}^0$ so that
$\text{div}(w) = p$ in $D$. Specifically, the function $w$ is
defined in the following way
\begin{equation*}
w = \sum_{i=1}^{N_0} w_i \Psi^{i,\text{off}}_{r_i}, \quad\quad w_i =
\int_{E_i} \eta\cdot m_i
\end{equation*}
and, in this proof only, we normalize the basis functions so that
$\int_{E_i} \Psi^{i,\text{off}}_{r_i} \cdot m_i = 1$. Thus,
\begin{equation}
\int_D p^2 = \int_D \text{div}(\eta) p = \sum_{i=1}^{N_0} \int_{E_i}
(\eta\cdot m_i) [p] = \sum_{i=1}^{N_0} \int_{E_i} w_i
(\Psi^{i,\text{off}}_{r_i}\cdot m_i) [p] = \int_D \text{div}(w) p,
\label{eq:infsup_a}
\end{equation}
where $[p]$ is the jump of $p$ across the coarse edge.

To show (\ref{eq:infsup1}), it remains to estimate $\|
w\|_{\kappa^{-1},D}$. Notice that,
\begin{equation*}
\| w\|_{\kappa^{-1},D}^2 = \int_D \kappa^{-1} w\cdot w \leq
\sum_{i=1}^{N_0} \int_{\omega_i} \kappa^{-1} w_i^2
\Psi^{i,\text{off}}_{r_i} \cdot \Psi^{i,\text{off}}_{r_i}.
%\leq \sum_{i=1}^{N_0} \frac{w_i^2}{\lambda^{(i)}_{r_i}} \int_{E_i} ( \Psi^{i,\text{off}}_{r_i} \cdot m_i)^2
\end{equation*}
For each $i$, we have $w_i^2 \leq H \int_{E_i} (\eta\cdot m_i)^2$.
Thus,
\begin{equation*}
\| w\|_{\kappa^{-1},D}^2 \preceq H \Big( \max_{1\leq i\leq N_0}
\int_{\omega_i} \kappa^{-1}  \Psi^{i,\text{off}}_{r_i} \cdot
\Psi^{i,\text{off}}_{r_i} \Big) \sum_{K\in\mathcal{T}^H}
\int_{\partial K}  (\eta\cdot n)^2.
\end{equation*}
Since the above inequality holds for any $\Psi^{i,\text{off}}_{r}$
such that $\int_{E_i} \Psi^{i,\text{off}}_r \cdot m_i \ne 0$, we have
\begin{equation}
\| w\|_{\kappa^{-1},D}^2 \preceq H \Big( \max_{1\leq i\leq N_0} \min_{r}
\int_{\omega_i} \kappa^{-1}  \Psi^{i,\text{off}}_{r} \cdot
\Psi^{i,\text{off}}_{r} \Big) \sum_{K\in\mathcal{T}^H}
\int_{\partial K}  (\eta\cdot n)^2,
\label{eq:infsup_b}
\end{equation}
where the above minimum is taken over all indices $r$ with the property
$\int_{E_i} \Psi^{i,\text{off}}_r \cdot m_i \ne 0$.

Finally, we will estimate $\int_{\partial K} (\eta\cdot n)^2$ for
every coarse grid block $K$. By the Green's identity, we have
\begin{equation*}
\int_{\partial K} (\eta\cdot n) z = \int_K \nabla \zeta\cdot \nabla
\tilde{z} + \int_K p \tilde{z},
\end{equation*}
where $z\in H^{\frac{1}{2}}(\partial K)$ and $\tilde{z}\in H^1(K)$
is any extension of $z$ in $K$. By Cauchy-Schwarz inequality,
\begin{equation*}
\begin{split}
\int_{\partial K} (\eta\cdot n) z
&= \int_K \nabla \zeta\cdot \nabla \tilde{z} + \int_K p \tilde{z} \\
&\preceq \Big( \| \nabla \zeta \|^2_{L^2(K)} + \| p \|_{L^2(K)}^2 \Big)^{\frac{1}{2}} \| \tilde{z}\|_{H^1(K)} \\
& \leq C_K \Big( \| \nabla \zeta \|^2_{L^2(K)} + \| p \|_{L^2(K)}^2
\Big)^{\frac{1}{2}} \| z\|_{H^{\frac{1}{2}}(\partial K)},
\end{split}
\end{equation*}
where the constant $C_K$ depends on $K$.
Thus,
\begin{equation*}
\int_{\partial K} (\eta\cdot n)^2 \leq C^2_K \Big( \| \nabla \zeta
\|^2_{L^2(K)} + \| p \|_{L^2(K)}^2 \Big).
\end{equation*}
By a scaling argument, we obtain
\begin{equation*}
H \int_{\partial K} (\eta\cdot n)^2 \preceq  \| \nabla \zeta
\|^2_{L^2(K)} + \| p \|_{L^2(K)}^2.
\end{equation*}
Summing the above over all coarse grid blocks $K$ and using $\| \nabla \zeta\|^2_{L^2(K)} \preceq \|p\|_{L^2(K)}^2$,
we have $H \sum_{K\in\mathcal{T}^H} \int_{\partial K} (\eta\cdot n)^2 \preceq \| p\|^2_{L^2(K)}$.
Hence, we obtain the desired bound (\ref{eq:infsup1})
by using (\ref{eq:infsup_a}) and (\ref{eq:infsup_b}).

\begin{flushright}
$\square$
\end{flushright}

Now we state and prove the convergence theorem for the mixed GMsFEM
(\ref{approximation-problem}).

\begin{theorem}
\label{thm1} Let $v_h$ be the fine-grid solution obtained in
(\ref{fine-problem}) and $v_H$ be the mixed GMsFEM solution obtained
in (\ref{approximation-problem}). Then, the following estimate holds
\begin{equation}
\label{estimate} \int_D \kappa^{-1} | v_h - v_H |^2 \preceq
C_{\text{\rm infsup}}^2 \Lambda^{-1} \sum_{i=1}^{N_0} a_i(\widehat{v},\widehat{v}) +
 \max_{K\in\mathcal{T}^H} \Big( \kappa^{-1}_{\text{min},K} \Big) \sum_{i=1}^{N_0} \| f - \overline{f}\|_{L^2(K_i)}^2,
\end{equation}
where $\Lambda = \min_{1\leq i\leq N_0} \lambda^{(i)}_{l_i+1}$ and
$\widehat{v}$ is the projection of $v_h$ defined in
(\ref{proj-v})-(\ref{proj-v1}).
\end{theorem}

{\it Proof}. Subtracting (\ref{fine-problem}) by
(\ref{approximation-problem}), and using the fact that
$V^0_{\text{off}}\subset V^0_h$ and $Q_H\subset Q_h$, we have
\begin{eqnarray}
\label{error-diff}
\begin{split}
\int_D \kappa^{-1} (v_h- v_H) \cdot w_H  - \int_D \text{div} (w_H)
(p_h-p_H) &=0,\quad
\forall w_H\in  {V}_{\text{off}}^0, \\
\int_D \text{div} (v_h-v_H) q_H  &=0, \quad \forall q_H \in Q_H.
\end{split}
\end{eqnarray}
By (\ref{error-eqn}), for each coarse-grid block $K$, we have
\begin{equation*}
\int_K \text{div}( v_h - \widehat{v} ) q_H = \int_K (f-\overline{f})
q_H = 0, \quad \forall q_H \in Q_H
\end{equation*}
since $q_H$ is a constant function on $K$. Similarly, since
$\text{div}(w_H)$ is a constant function for any $w_H\in
V^0_{\text{off}}$, by (\ref{proj-v1}), we have
\begin{equation*}
\int_D \text{div} (w_H) p_h = \int_D \text{div} (w_H) \widehat{p}.
\end{equation*}
Thus, (\ref{error-diff}) can be written as
\begin{eqnarray}
\label{error-diff-1}
\begin{split}
\int_D \kappa^{-1} (v_h- v_H) \cdot w_H  - \int_D \text{div} (w_H)
(\widehat{p}-p_H) &=0,\quad
\forall w_H\in  {V}_{\text{off}}^0, \\
\int_D \text{div} (\widehat{v}-v_H) q_H  &=0, \quad \forall q_H \in
Q_H.
\end{split}
\end{eqnarray}
Notice that $\widehat{v}\in V_{\text{snap}}$. We can therefore write
$\widehat{v}$ as
\begin{equation}
\widehat{v} = \sum_{i=1}^{N_0} \sum_{k=1}^{J_i} \widehat{v}_{ij}
\Psi^{i,\text{off}}_{k}. \label{eq:v1}
\end{equation}
We then define $\widehat{v}_{\text{off}} \in V_{\text{off}}$ by
\begin{equation}
\widehat{v}_{\text{off}} = \sum_{i=1}^{N_0} \sum_{k=1}^{l_i}
\widehat{v}_{ij} \Psi^{i,\text{off}}_{k}, \label{eq:v2}
\end{equation}
where we recall that $l_i \leq J_i$ is the number of eigenfunctions
selected for the coarse neighborhood $\omega_i$. Notice that
%$\widehat{v}_{ij}=0$ when the coarse edge $E_i$ belongs to the
%global boundary $\partial D$. Thus,
$\widehat{v}_{\text{off}} \in
V_{\text{off}}^0$. We can further write (\ref{error-diff-1}) as
\begin{eqnarray}
\label{error-diff-2}
\begin{split}
\int_D \kappa^{-1} (v_h- v_H) \cdot w_H  - \int_D \text{div} (w_H)
(\widehat{p}-p_H) &=0,\quad
\forall w_H\in  {V}_{\text{off}}^0, \\
\int_D \text{div} (\widehat{v}_{\text{off}}-v_H) q_H  &=\int_D
\text{div} (\widehat{v}_{\text{off}}-\widehat{v}) q_H, \quad \forall
q_H \in Q_H.
\end{split}
\end{eqnarray}
Taking $w_H = \widehat{v}_{\text{off}}-v_H$ and $q_H =
\widehat{p}-p_H$ in (\ref{error-diff-2}), and adding the resulting
equations, we obtain
\begin{equation}
\label{error-diff-3}
\int_D \kappa^{-1} (v_h- v_H) \cdot (\widehat{v}_{\text{off}}-v_H)
 =\int_D \text{div} (\widehat{v}_{\text{off}}-\widehat{v}) (\widehat{p}-p_H)
\end{equation}
%Using (\ref{proj-v1}) and the fact that $\text{div}
%(\widehat{v}_{\text{off}}-\widehat{v})$ is a constant function on
%each coarse-grid block, the above equation becomes
%\begin{equation}
%\label{error-diff-3} \int_D \kappa^{-1} (v_h- v_H) \cdot
%(\widehat{v}_{\text{off}}-v_H)
 %=0.
%\end{equation}
By the inf-sup condition (\ref{eq:infsup1}) and the error equation (\ref{error-diff-2}),
we have $$\| \widehat{p}-p_H\|_{L^2(D)} \preceq C_{\text{infsup}} \| v_h - v_H \|_{\kappa^{-1},D}.$$
Moreover, by the definition of the spectral problem (\ref{eq:pod1}), we have
\begin{equation*}
\int_D (\text{div}( \widehat{v}_{\text{off}} - \widehat{v}) )^2
\preceq \sum_{i=1}^{N_0} \int_{\omega_i} (\text{div}( \widehat{v}_{\text{off}} - \widehat{v}) )^2
\preceq \sum_{i=1}^{N_0} s_i( \widehat{v}_{\text{off}} - \widehat{v}, \widehat{v}_{\text{off}} - \widehat{v} ).
\end{equation*}
We can then derive from (\ref{error-diff-3}) the following
\begin{equation*}
\| v_h - v_H \|_{\kappa^{-1},D}^2 \preceq \| \widehat{v}_{\text{off}}- v_h \|^2_{\kappa^{-1},D}
+ C_{\text{infsup}}^2 \sum_{i=1}^{N_0} s_i( \widehat{v}_{\text{off}} - \widehat{v}, \widehat{v}_{\text{off}} - \widehat{v} ).
\end{equation*}
%which implies
%\begin{equation*}
%\label{error-diff-4}
%\int_D \kappa^{-1} (\widehat{v}_{\text{off}}-v_H) \cdot
%(\widehat{v}_{\text{off}}-v_H)
 %=\int_D \kappa^{-1} (\widehat{v}_{\text{off}}-v_h) \cdot
 %(\widehat{v}_{\text{off}}-v_H) + \int_D \text{div} (\widehat{v}_{\text{off}}-\widehat{v}) (\widehat{p}-p_H).
%\end{equation*}
%Therefore, we obtain
%\begin{equation*}
%\| \widehat{v}_{\text{off}}-v_H \|_{\kappa^{-1},D} \preceq \|
%\widehat{v}_{\text{off}}-\widehat{v} \|_{\kappa^{-1},D} + \|
%\widehat{v}-v_h \|_{\kappa^{-1},D}
%\end{equation*}
Using the triangle inequality $\| \widehat{v}_{\text{off}} - v_h \|_{\kappa^{-1},D}
\leq \| \widehat{v}_{\text{off}} - \widehat{v} \|_{\kappa^{-1},D} + \| \widehat{v}-v_h\|_{\kappa^{-1},D}$
and
%To estimate the term $\| \widehat{v}_{\text{off}}-\widehat{v}
%\|_{\kappa^{-1},D}$, we notice that
\begin{equation*}
\| \widehat{v}_{\text{off}}-\widehat{v} \|_{\kappa^{-1},D}^2 \preceq
\sum_{i=1}^{N_0} \| \widehat{v}_{\text{off}}-\widehat{v}
\|_{\kappa^{-1},\omega_i}^2 \preceq \sum_{i=1}^{N_0}
s_i(\widehat{v}_{\text{off}}-\widehat{v},
\widehat{v}_{\text{off}}-\widehat{v}),
\end{equation*}
we obtain
\begin{equation}
\| v_h - v_H \|_{\kappa^{-1},D}^2 \preceq \| \widehat{v}- v_h \|^2_{\kappa^{-1},D}
+ C_{\text{infsup}}^2 \sum_{i=1}^{N_0} s_i( \widehat{v}_{\text{off}} - \widehat{v}, \widehat{v}_{\text{off}} - \widehat{v} ).
\label{error-diff-4}
\end{equation}
The first term on the right hand side of (\ref{error-diff-4}) can be estimated
by Lemma \ref{lemma1}.
For the second term on the right hand side of (\ref{error-diff-4}),
by (\ref{eq:v1})-(\ref{eq:v2}) and the fact that
$\Psi_k^{i,\text{off}}$ are eigenfunctions of (\ref{eq:localspec}),
we have
\begin{equation*}
s_i(\widehat{v}_{\text{off}}-\widehat{v},
\widehat{v}_{\text{off}}-\widehat{v}) = \sum_{k=l_i+1}^{J_i}
(\lambda_{k}^{(i)})^{-1} (\widehat{v}_{ik})^2 a_i(
\Psi_k^{i,\text{off}}, \Psi_k^{i,\text{off}}).
\end{equation*}
By the ordering of the eigenvalues (\ref{eq:eigenorder}) and
orthogonality of eigenfunctions, we obtain
\begin{equation*}
s_i(\widehat{v}_{\text{off}}-\widehat{v},
\widehat{v}_{\text{off}}-\widehat{v}) \leq
(\lambda_{l_i+1}^{(i)})^{-1} a_i( \widehat{v}_{\text{off}}
-\widehat{v}, \widehat{v}_{\text{off}}-\widehat{v}) \leq
(\lambda_{l_i+1}^{(i)})^{-1} a_i( \widehat{v}, \widehat{v}).
\end{equation*}
Combining the above results, we have
\begin{equation*}
\sum_{i=1}^{N_0} s_i( \widehat{v}_{\text{off}} - \widehat{v}, \widehat{v}_{\text{off}} - \widehat{v} ) \leq
\sum_{i=1}^{N_0} (\lambda_{l_i+1}^{(i)})^{-1} a_i( \widehat{v},
\widehat{v}).
\end{equation*}
This completes the proof.

%Finally, the desired estimate can be obtained by using the following
%triangle inequality
%\begin{equation*}
%\| v_h - v_H \|_{\kappa^{-1},D} \leq \|v_h -
%\widehat{v}\|_{\kappa^{-1},D} + \| \widehat{v} -
%\widehat{v}_{\text{off}} \|_{\kappa^{-1},D} + \|
%\widehat{v}_{\text{off}} - v_H \|_{\kappa^{-1},D},
%\end{equation*}
%with the above bound for $\| \widehat{v}_{\text{off}}-\widehat{v}
%\|_{\kappa^{-1},D}$ and Lemma \ref{lemma1}.

\begin{flushright}
$\square$
\end{flushright}

We remark that in the error estimate (\ref{estimate}), the first and
second terms on the right-hand-side represent the errors due to the
spectral basis functions and the coarse grid discretization,
respectively.

%\subsubsection{Error analysis for spectral problem 2}
%Let $E_i$ be coarse edge with corresponding neighborhood $\omega_i$.
%We consider
%\begin{equation}
%v \cdot n = \lambda [\phi] \quad\quad \text{ on} \quad E_i
%\end{equation}
%only on the divergence free subspace.
%The notation $[\phi]$ is the jump of $\phi$.
%This means we always include the constant function (on edge) as basis.
%Then we have
%\begin{equation}
%\begin{split}
%&\: \int_{\omega_i} \kappa^{-1} (\widehat{v}- v_H)\cdot (\widehat{v}-v_H) + \int_{\omega_j} \kappa^{-1} (\widehat{v}- v_H)\cdot (\widehat{v}-v_H) \\
%= &\: \int_{\omega_i} (\widehat{v}- v_H)\cdot  \nabla (\widehat{\phi}-\phi_H) + \int_{\omega_j} (\widehat{v}- v_H)\cdot  \nabla (\widehat{\phi}-\phi_H) \\
%= &\: \int_{E} (\widehat{v}- v_H)\cdot n \, [ \widehat{\phi}-\phi_H ] \\
%\preceq &\: \frac{1}{\lambda_{i,L_i+1}} \int_{E} ((\widehat{v}- v_H)\cdot n )^2 \\
%\preceq &\: \frac{1}{\lambda_{i,L_i+1}} \int_{E} (\widehat{v}\cdot n )^2
%\end{split}
%\end{equation}
%
%Numerically, the above eigenvalue problem can be implemented as
%\begin{equation}
%\int_{\omega_i} \kappa^{-1} v \cdot u = \lambda \int_{E_i} [\phi_v]
%[\phi_u]
%\end{equation}

%\section{Relation to homogenization}

%- Eigenvalue structure and the discussion on the first basis function

%- Relation to homogenization and MsFEM

\section{Oversampling approach}
\label{sec:ovs}

One can use an oversampling approach to improve the accuracy of the
method. The main idea of the oversampling method is to use larger
domains to compute snapshots. Furthermore, performing POD in the
snapshot space, we can achieve a lower dimensional approximation
space.
Oversampling technique can be particularly helpful for problems
with scale separation. This is because by
taking the restriction of the local solutions in larger domains
in the interior, we avoid the pollution
effects near the boundaries.

Let $\Omega$ be a conforming subset of $D$. By conforming subset, we
mean that $\Omega$ is formed by the union of connected fine grid
elements. For a given function $\psi$ defined on $\partial\Omega$,
we denote $(\mathcal{H}_{\Omega}(\psi), \phi_{\Omega}) \in
V_h(\Omega) \times Q_h(\Omega)$ by the solution of the weak form of
the following problem
\begin{eqnarray}
\label{harmonic}
\begin{split}
 \kappa^{-1} \mathcal{H}_{\Omega}(\psi) +\nabla
\phi_{\Omega} &=0 \quad \text{in} \quad \Omega, \\
\text{div} ( \mathcal{H}_{\Omega}(\psi) ) &= c_{\Omega} \quad \text{in} \quad \Omega, \\
\mathcal{H}_{\Omega}(\psi) \cdot n &= \psi \quad \text{on} \quad
\partial\Omega,
\end{split}
\end{eqnarray}
where $c_{\Omega} = |\Omega|^{-1} \int_{\partial \Omega} \psi$,
$V_h(\Omega)$ and $Q_h(\Omega)$ are the restrictions of $V_h$ and
$Q_h$ on $\Omega$ respectively. We call $\mathcal{H}_{\Omega}(\psi)$
the $\kappa$-harmonic extension of $\psi$ in $\Omega$.

Let $E_i \in \mathcal{E}^H$ be an interior coarse edge, and let
$\omega_i^{+}$ be a conforming subset of $D$ with $E_i$ lying in the
interior of $\omega_i^{+}$, see {\sc Fig}.~\ref{fig:coarse_grid} for an
example of $\omega_i^{+}$. Let $W_i(\partial \omega_i^{+})$ be the
set of all piecewise constant functions defined on
$\partial\omega_i^{+}$ with respect to the fine grid partition.
Consider the following set of functions defined on $E_i$
\begin{equation*}
\Big\{ \mathcal{H}_{\omega_i^{+}}(\psi_j) \cdot m_i |_{E_i} \quad ,
\quad \psi_j \in W_i(\partial\omega_i^{+}) \Big\}.
\end{equation*}
By performing a standard POD on the above space, and selecting the
first $l_i^{+}$ dominant modes $\psi_j^{i,\text{ovs}}$, we obtain
the following space
\begin{equation*}
V_{\text{ovs}}(E_i) = \text{span} \Big\{ \Psi_j^{i,\text{ovs}} \; ,
\; 1\leq j \leq l_i^{+} \Big\},
\end{equation*}
where the basis functions $\Psi_j^{i,\text{ovs}}$ are obtained by
solving (\ref{MMs-basis-eq}) with the boundary condition
(\ref{eq:BC}) replaced by $\Psi_j^{i,\text{ovs}} \cdot m_i =
\psi_j^{i,\text{ovs}}$ on $E_i$. We call this local oversampling
space. The oversampling space $V_{\text{ovs}}$ is obtained by the
linear span of all local oversampling spaces. To obtain a numerical
solution, we solve (\ref{approximation-problem}) with $V_H =
V_{\text{ovs}}$.

Next, we discuss the outline of the convergence analysis
for the oversampling approach. For
any $v_h\in V_h$ and for every $E_i\in\mathcal{E}^H$, we define
\eric{$\mu_{E_i}$} as
\[
\eric{\mu_{E_i}}=\mathcal{H}_{\omega_i^{+}}(v_h\cdot n|_{\partial
\omega_i^{+}} )\cdot \eric{m_i}|_{E_i},
\]
\eric{which is the normal component on $E_i$ of the $\kappa$-harmonic extension of $v_h\cdot n$
in the oversampled region $\omega_i^{+}$.
Using $\mu_{E_i}$, we can then define $\widetilde{v}\in V_{\text{snap}}$ by}
\begin{equation*}
\eric{\widetilde{v} = \sum_{i=1}^{N_e}  \mu_{E_i} \Psi^{(i)}_0},
\end{equation*}
\eric{where $\Psi^{(i)}_0 \in V_{\text{snap}}^{(i)}$ satisfies $\Psi^{(i)}_0 \cdot m_i = 1$ on $E_i$.}
%\remove{Denote $\widetilde{v}_{E_i}=\mu_{E_i}$ over all edges by
%$\widetilde{v}_\Gamma$ (normal trace), $\widetilde{v}_\Gamma=\sum_i
%\widetilde{v}_{E_i} \chi_{E_i}$, where $\chi_{E_i}$ is $1$ on edge
%$E_i$ and 0 all other edges. We consider
%$\widetilde{v}=\mathcal{H}(\widetilde{v}_\Gamma)=\sum_{E_i,K\bigcap
%E_i\not = \varnothing}
%\mathcal{H}_K(\widetilde{v}_{E_i}\chi_{E_i})$.
%First,
%we will estimate
%\[
%\|\widehat{v}-\widetilde{v}\|_{0,D},
%\]
%}
%Note $\widehat{v}=\sum_{E_i,K\bigcap E_i\not = \varnothing}
%\mathcal{H}_K(v_{E_i}\chi_{E_i})$, where $v_{E_i}=v_h\cdot n|_{E_i}$.
 Next, we have
\begin{equation}
\begin{split}
\|\widehat{v}-\widetilde{v}\|_{\kappa^{-1},D}^2 &=
%\sum_K \|\sum_{E_i,K\bigcap E_i\not = \varnothing} \mathcal{H}_K(v_{E_i}\chi_{E_i}-\mu_{E_i}\chi_{E_i})\|_{\kappa^{-1},K}^2\preceq\\
%\sum_K \|\sum_{E_i,K\bigcap E_i\not = \varnothing} v_{E_i}\Psi_0^{(i)}-\mu_{E_i}\Psi_0^{(i)}\|_{\kappa^{-1},K}^2\preceq\\
\sum_{i=1}^{N_e} \| (\widehat{v}\cdot m_i) \Psi^{(i)}_0 - \mu_{E_i} \Psi^{(i)}_0 \|_{\kappa^{-1},\omega_i}^2 \\
&\preceq \sum_{i=1}^{N_e} H\| (\widehat{v}\cdot m_i) - \mu_{E_i} \|_{L^2(E_i)}^2\preceq
H\delta,
\end{split}
\end{equation}
%where $v_{E_i}$ is the normal trace of $v_h$.
where we assumed that $H\| (\widehat{v}\cdot m_i) - \mu_{E_i} \|_{L^2(E_i)}^2\preceq H\delta$ and $\widehat{v}\cdot m_i$ is the normal trace
of $\widehat{v}$ on $E_i$.
%\[
%\remove{\|v_{E_i}-\widetilde{v}_{E_i}\|_{0,E_i}^2\preceq H\delta.}
%\]
\eric{If the forcing is constant within the union of $\omega_i^{+}$ and $\omega_i$, then $\delta=0$.
Otherwise, this value depends on the smoothness of $\kappa$ and $f$.
For homogenization problems, one can show that $\delta$ is small.}

Next, we choose an
appropriate interpolant \eric{$\widetilde{v}_{\text{ovs}}$} and compare it with $\widetilde{v}$.
Note that, we can write
\begin{equation*}
\widetilde{v}_{\text{ovs}} = \sum_{i=1}^{N_e} \sum_{j=1}^{l_i^+}c_j^i\Psi_j^{i,\text{ovs}}
\end{equation*}
for some constants $c_j^i$. Therefore,
\begin{equation}
\begin{split}
\|\widetilde{v}-\widetilde{v}_{\text{ovs}}\|_{\kappa^{-1},D}^2 &=\|\sum_{i=1}^{N_e}
\Big( \mu_{E_i} \Psi^{(i)}_0 - \sum_{j=1}^{l_i^+}c_j^i\Psi_j^{i,\text{ovs}} \Big)\|_{\kappa^{-1},D}^2  \\
&\preceq \sum_{i=1}^{N_e}
H\|\mu_{E_i}-\sum_{j=1}^{l_i^+}c_j^i\Psi_j^{i,\text{ovs}}\cdot
m_i\|_{L^2(E_i)}^2.
%\preceq\\
%  \sum_{i=1}^{N^e} {1\over \lambda^+_{l_i^++1}} H\|\widetilde{v}_{E_i}-\sum_{j=1}^{l_i^+}c_j^i\Psi_j^{i,\text{over}}\cdot n\|_{0,\partial \omega_i^{+}}^2\preceq {1\over \Lambda} \sum_{i=1}^{N^e} H \|\widetilde{v}_{E_i}\|_{0,\partial \omega_i^{+}}^2.
\end{split}
\end{equation}
Denote by $\Psi^E$, the restriction of the snapshots on the edge
$E$.
% $\mathcal{H}_{\omega_i^+}(\psi_j)\cdot m_i$ restricted on $E_i$.
We would like to find a reduced dimensional representation of
$\Psi^E$ such that $\|\Psi^E - \Phi^E_r C_r\|$ is small, where
$\Phi^E_r$ is the reduced-dimensional representation (the matrix of
the size $N^e\times N^r$), where $N^r$ is the reduced dimensional
and $C^r$ is the matrix of the size $N^r\times N^{\partial
\omega_i^+}$,
where $N^{\partial
\omega_i^+}$ is the number of fine-grid edges on $\partial
\omega_i^+$.
 This is achieved by POD as described above and we have
 $\|\Psi^E  - \Phi^E_r C_r \|_F\leq 1/\lambda_{l_i^+ + 1}^+ $. From here,
one can show that given values of the velocity $z$ on the boundary
of $\partial \omega_i^+$, we have $\|\Psi^E z - \Phi^E_r C_r
z\|_2\leq (1/\lambda_{l_i^+ + 1}^+) \|z\|_2$. Combining these
estimates, we have
\begin{equation}
\begin{split}
\|\widetilde{v}- \widetilde{v}_{\text{ovs}} \|_{\kappa^{-1},D}^2
%=\|\sum_{i=1}^{N^e}  \mathcal{H}_K(\widetilde{v}_{E_i}-\sum_{j=1}^{l_i^+}c_j^i\Psi_j^{i,\text{over}}\cdot m_i)\|_{0,D}^2\preceq \sum_{i=1}^{N^e} H\|\widetilde{v}_{E_i}-\sum_{j=1}^{l_i^+}c_j^i\Psi_j^{i,\text{over}}\cdot m_i\|_{0,E_i}^2
\preceq
 % \sum_{i=1}^{N^e} {1\over \lambda^+_{l_i^++1}} H\|\widetilde{v}_{E_i}-\sum_{j=1}^{l_i^+}c_j^i\Psi_j^{i,\text{over}}\cdot n\|_{0,\partial \omega_i^{+}}^2\preceq {1\over \Lambda} \sum_{i=1}^{N^e} H \|\widetilde{v}_{E_i}\|_{0,\partial \omega_i^{+}}^2.
 \sum_{i=1}^{N_e} {1\over \lambda^+_{l_i^++1}} H\|v_h\|_{L^2)\partial \omega_i^{+})}^2\preceq {1\over \Lambda^+} \sum_{i=1}^{N_e} H \|v_h\|_{L^2(\partial \omega_i^{+})}^2
\end{split}
\end{equation}
where $\Lambda^+ = \min \{ \lambda_{l_i^+ +1}^+ \}$.
%In the last step, we have used the fact that
%$\sum_{j=1}^{l_i^+}c_j^i\Psi_j^{i,\text{over}}$ is chosen to be the
%interpolant of $\widetilde{v}$ in $V_H^{\text{over}}$. \\
%(\textbf{need some estimate in first step to go from $v_H$ to
%$\sum_{j=1}^{l_i^+}c_j^i\Psi_j^{i,\text{over}}$})

One can consider an alternative approach where the snapshot space is
obtained by performing POD as described above. More precisely, we
use $V_{\text{ovs}}$ as the snapshot space that can have a lower
dimension compared to the original snapshot space that corresponds
to non-oversampling case. As a next step, we perform a spectral
decomposition following the non-oversampling case by considering
$\mathcal{H}_K(\Psi_j^{i,\text{ovs}})$ as a snapshot space. We denote
this snapshot space by  $V_{\text{ovs}}^R$, where $R$ stands for
reduced dimension.
 The main advantage of this approach is that
a lower dimensional snapshot space is used in the spectral
decomposition and this snapshot space allows achieving a low
dimensional structure when the problem has a scale separation. The
latter may not hold if we apply non-oversampling procedure.
%In this approach, we first truncate the snapshot space based on POD and then apply our local spectral problem
%to define multiscale basis functions. The main advantage of this approach is that the snapshot space
%can be drastically reduced for some cases, e.g., the case with scale separation.
To obtain the convergence analysis,
%we consider the snapshot space obtained by selecting dominant modes
%in the snapshot space. In particular,
we  show that for every $v_h \in V_h$, there exists
$\mu_{E_i}^{R} \Psi^{(i)}_0 \in V_{\text{ovs}}^R $ in the snapshot
space, such that
\begin{equation}
\label{eq:reduced_snapshot} \|\mu_{E_i} -
\mu_{E_i}^R \|_{L^2(E_i)}^2\preceq 1/\lambda_{l_i^+ + 1}^+,
\end{equation}
where $1/\lambda_{l_i^+ + 1}^+$ is the lowest eigenvalue that the
corresponding eigenvector is not included in the snapshot space.
This follows from standard POD result which provides an estimate for
$\|\Psi^E  - \Phi^E_r C_r \|_F\leq 1/\lambda_{l_i^+ + 1}^+$. Under
this condition and using the fact that $\|A z \|_2\leq \|A\|_F
\|z\|_2$, we obtain (\ref{eq:reduced_snapshot}). Using this reduced
snapshot space, we can repeat our previous argument in Section
\ref{sec:analysis} and obtain the convergence rate.

\section{Numerical results}
\label{sec:numresults}

\begin{figure}[htp!]
\centering \subfigure[$\kappa_1$]{
\includegraphics[scale=.21]{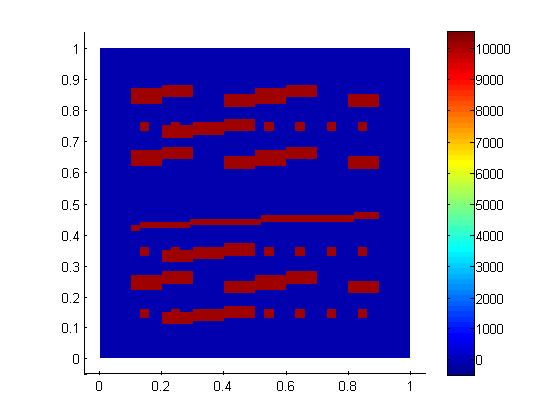}}
\subfigure[$\kappa_2$]{
\includegraphics[scale=.21]{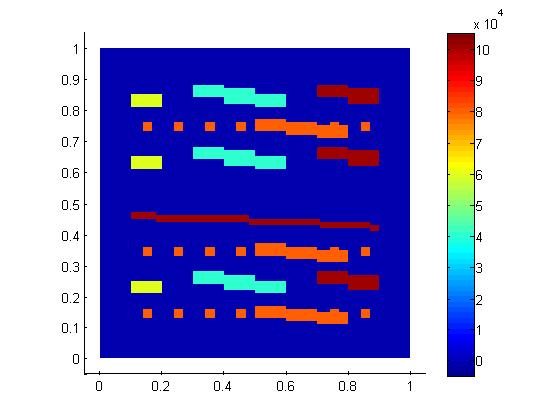}}
\subfigure[$\kappa_3$ in log$_{10}$ scale]{
\includegraphics[scale=.21]{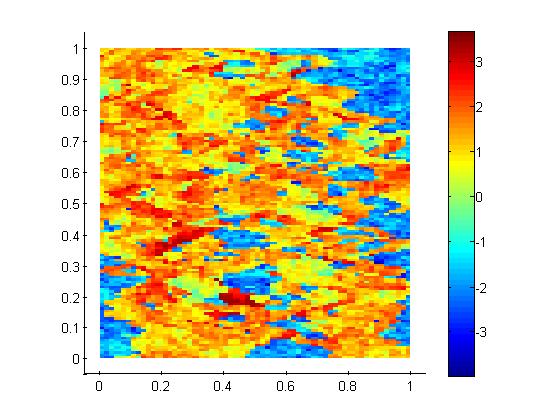}}
  \caption{Three permeability fields in the numerical experiments}\label{fig:coef}
\end{figure}

In this section, we will present some numerical results to show the
performance of the mixed GMsFEM (\ref{approximation-problem}) for
approximating the flow problem (\ref{pv-system}). In all simulations
reported below, the computational domain $D = (0,1)^2$.
%throughout this section, and we solve problem (\ref{pv-system})
%using our mixed generalized multiscale finite element method
%(\ref{approximation-problem}).
%In order to see the robustness of the
%method,
The coarse grid $\mathcal{T}^H$ and the fine grid $\mathcal{T}^h$
are $N\times N$ and
 $n\times n$ uniform meshes, respectively.
 A fixed fine-grid size with $n=200$ is employed.
Moreover, we will consider three different permeability fields
$\kappa$, as depicted in {\sc Fig}.~\ref{fig:coef}.
These permeability fields have the same resolution as the fine-grid
size.
%All the grids in the numerical experiments are uniform rectangular
%grids. We partition the domain into a $n$ by $n$ fine grid and $N$
%by $N$ coarse grid.
We will present the performance of the mixed GMsFEM for three types
of applications; namely, we present single-phase flow problems,
single-phase flow and transport problems, and two-phase flow and
transport problems. To facilitate the presentation, we let
$(v_f,p_f), (v_s,p_s)$ and $(v_o,p_o)$ be the fine-grid solution,
snapshot solution, and the GMsFEM solution respectively, where the
snapshot solution is the solution of the discrete system
(\ref{approximation-problem}) with all basis functions in the
snapshot space are selected. Notice that the snapshot solution
contains only the coarse-grid discretization error and the GMsFEM
solution contains both coarse-grid and spectral errors, see Theorem
\ref{thm1}. Furthermore, we define the following error quantities
for the velocity field
\begin{equation*}
E_{of}(v) := \|v_o-v_f\|_{\kappa^{-1},D} / \|v_f\|_{\kappa^{-1},D},
\;\;  E_{os}(v) :=
\|v_o-v_s\|_{\kappa^{-1},D}/\|v_s\|_{\kappa^{-1},D},
\end{equation*}
which we term the total error and the spectral error, respectively.
For pressure, we define the corresponding error quantities by
\begin{equation*}
E_{of}(p) := \|p_o-p_f\|_{L^2(D)} / \|p_f\|_{L^2(D)}, \;\; E_{os}(p)
:= \|p_o-p_s\|_{L^2(D)}/\|p_s\|_{L^2(D)}.
\end{equation*}
These error quantities are used to measure the performance of the
mixed GMsFEM in the examples below.

\subsection{Single-phase flow}

We consider single-phase flow in this section. For the simulations,
we will use two different coarse-mesh sizes with $N=10$ and $N=20$,
called case $1$ and case $2$, respectively. The numerical results
for the permeability fields $\kappa_1$ and $\kappa_2$, as well as
the use of the above two spectral problems (\ref{eq:pod1}) and
(\ref{eq:pod2}) are shown in {\sc Tables} \ref{TB1}-\ref{TB4}. In
these tables, the term "dof per $E$" means the number of basis
functions used for that coarse edge $E$. We remark that, for
spectral problem 2, the first eigenfunction is always taken as the
field with constant normal component on $E_i$.
%We denote the fine scale solution, snapshot
%solution and the offline solution by $(v_f,p_f), (v_s,p_s),
%(v_o,p_o)$, respectively. The weighted $L^2$ error of $v_o$ relative
%to $v_f$ is computed as $\|v_o-v_f\|_\kappa/\|v_f\|_\kappa$, where
%$\|v\|_\kappa = (\int_\Omega \kappa^{-1} v^2 dx)^{1/2}$. Similarly,
%the weighted $L^2$ error of $v_o$ relative to $v_s$ is computed as
%$\|v_o-v_s\|_\kappa/\|v_s\|_\kappa$. In order to put the tables in a
%more compact way, we adopt the following notation:
%\begin{equation*}
%E_{of}(v) := \|v_o-v_f\|_\kappa/\|v_f\|_\kappa, \;\;  E_{os}(v) :=
%\|v_o-v_s\|_\kappa/\|v_s\|_\kappa.
%\end{equation*}
%$E_{of}(p)$ and $E_{os}(p)$ are similarly defined.
In {\sc Tables} \ref{TB1}-\ref{TB2}, the convergence behaviors of
the method for the permeability field $\kappa_1$ are shown for cases
$1$ and $2$, respectively. Notice that, cases 1 and 2 decompose each
coarse-grid block as 20x20 and 10x10 grids, respectively. Therefore,
for each coarse edge, there are $20$ and $10$ basis functions for
cases $1$ and $2$ respectively. From these tables, we see clearly
the convergence of the method when basis functions are added to the
offline space. In addition, we see that the spectral errors
$E_{os}(v)$ and $E_{os}(p)$ converge to machine precision. On the
other hand, the total errors $E_{of}(v)$ and $E_{of}(p)$ converge to
a fixed error when the number of basis functions are increased. This
fixed error corresponds to the coarse grid discretization error and
cannot be improved by introducing more spectral basis functions.
Nevertheless, the coarse-grid error can be reduced by using a
smaller coarse mesh size. This is confirmed numerically in {\sc
Tables} \ref{TB1}-\ref{TB2}. In particular, when $N=10$, the level
of the coarse-grid error in velocity is about $2\%$; and when
$N=20$, the level of the coarse-grid error in velocity is reduced to
about $0.5\%$. We also observe a similar situation for pressure.
Regarding the results for the permeability field $\kappa_2$, the
results in {\sc Tables} \ref{TB3}-\ref{TB4} give a similar
conclusion.

%Two kinds of spectral problems are considered, they are introduced
%in section \ref{sec:offline}. For the ease of comparison of the
%performance between different spectral problems, in each of the
%convergence table below we fix the grid setting and coefficient, and
%put the errors corresponding to different spectral problem in one
%table.

\begin{table}[h!]
\footnotesize \caption{Convergence of the offline solution, $\kappa
= \kappa_1$, $n = 200$ and $N = 10$}\label{TB1} \centering
\begin{tabular}{|c||c|c|c|c|c|c|c|c|} \hline
& \multicolumn{4}{c|}{Spectral problem 1} & \multicolumn{4}{c|}{Spectral problem 2} \\
\hline dof per $E$ & $E_{of}(v)$ & $E_{of}(p)$ & $E_{os}(v)$ &
$E_{os}(p)$ & $E_{of}(v)$ & $E_{of}(p)$ & $E_{os}(v)$ & $E_{os}(p)$
\rule{0pt}{9pt}
\\  \hline
1 & 0.1331 & 0.0903 & 0.1329 & 0.0196 & 0.1523 & 0.1018 & 0.1525 & 0.0519\\
3 & 0.0569 & 0.0896 & 0.0535 & 0.0031 & 0.0840 & 0.0902 & 0.0823 & 0.0133 \\
5 & 0.0308 & 0.0898 & 0.0229 & 5.78e-04 & 0.0391 & 0.0898 & 0.0334 & 0.0031 \\
7 & 0.0236 & 0.0898 & 0.0112 & 1.39e-04 & 0.0278 & 0.0898 & 0.0186 & 0.0010 \\
9 & 0.0210 & 0.0898 & 0.0026 & 7.18e-06 & 0.0234 & 0.0898 & 0.0108 & 1.20e-04 \\
11 & 0.0208 & 0.0898 & 9.53e-13 & 4.87e-15 & 0.0208 & 0.0898 & 3.92e-13 & 4.94e-15 \\
20 & 0.0208 & 0.0898 & 3.92e-13 & 6.18e-15 & 0.0208 & 0.0898 & 3.96e-13 & 5.08e-15 \\
\hline
\end{tabular}
\end{table}

\begin{table}[h!]
\footnotesize \caption{Convergence of the offline solution, $\kappa
= \kappa_1$, $n = 200$ and $N = 20$}\label{TB2} \centering
\begin{tabular}{|c||c|c|c|c|c|c|c|c|} \hline
& \multicolumn{4}{c|}{Spectral problem 1} & \multicolumn{4}{c|}{Spectral problem 2} \\
\hline dof per $E$ & $E_{of}(v)$ & $E_{of}(p)$ & $E_{os}(v)$ &
$E_{os}(p)$ & $E_{of}(v)$ & $E_{of}(p)$ & $E_{os}(v)$ & $E_{os}(p)$
\rule{0pt}{9pt}
\\  \hline
1 & 0.1788 & 0.0601 & 0.1792 & 0.0373 & 0.1551 & 0.0677 & 0.1554 & 0.0483 \\
2 & 0.0460 & 0.0486 & 0.0459 & 0.0023 & 0.0861 & 0.0507 & 0.0861 & 0.0155 \\
3 & 0.0251 & 0.0486 & 0.0246 & 6.68e-04 & 0.0493 & 0.0488 & 0.0491 & 0.0055 \\
4 & 0.0115 & 0.0486 & 0.0102 & 1.15e-04 & 0.0233 & 0.0486 & 0.0227 & 0.0016 \\
5 & 0.0054 & 0.0486 & 3.47e-12 & 1.10e-14 & 0.0054 & 0.0486 & 4.29e-12 & 9.53e-15 \\
10 & 0.0054 & 0.0486 & 1.56e-12 & 1.29e-14 & 0.0054 & 0.0486 & 4.82e-13 & 9.61e-15 \\
\hline
\end{tabular}
\end{table}

\begin{table}[h!]
\footnotesize \caption{Convergence of the offline solution, $\kappa
= \kappa_2$, $n = 200$ and $N = 10$}\label{TB3} \centering
\begin{tabular}{|c||c|c|c|c|c|c|c|c|} \hline
& \multicolumn{4}{c|}{Spectral problem 1} & \multicolumn{4}{c|}{Spectral problem 2} \\
\hline dof per $E$ & $E_{of}(v)$ & $E_{of}(p)$ & $E_{os}(v)$ &
$E_{os}(p)$ & $E_{of}(v)$ & $E_{of}(p)$ & $E_{os}(v)$ & $E_{os}(p)$
\rule{0pt}{9pt}
\\  \hline
1 & 0.1404 & 0.0905 & 0.1403 & 0.0219 & 0.1482 & 0.0966 & 0.1482 & 0.0404\\
3 & 0.0561 & 0.0894 & 0.0526 & 0.0030 & 0.0778 & 0.0900 & 0.0757 & 0.0121 \\
5 & 0.0266 & 0.0896 & 0.0168 & 3.04e-04 & 0.0393 & 0.0897 & 0.0337 & 0.0047 \\
7 & 0.0232 & 0.0896 & 0.0105 & 1.20e-04 & 0.0277 & 0.0896 & 0.0185 & 0.0017 \\
9 & 0.0209 & 0.0896 & 0.0022 & 5.35e-06 & 0.0239 & 0.0896 & 0.0119 & 1.50e-04 \\
11 & 0.0208 & 0.0896 & 8.35e-13 & 8.19e-15 & 0.0208 & 0.0896 & 2.46e-11 & 7.48e-15 \\
20 & 0.0208 & 0.0896 & 4.98e-13 & 9.31e-15 & 0.0208 & 0.0896 & 5.00e-13 & 8.29e-15 \\
\hline
\end{tabular}
\end{table}

\begin{table}[h!]
\footnotesize \caption{Convergence of the offline solution, $\kappa
= \kappa_2$, $n = 200$ and $N = 20$}\label{TB4} \centering
\begin{tabular}{|c||c|c|c|c|c|c|c|c|} \hline
& \multicolumn{4}{c|}{Spectral problem 1} & \multicolumn{4}{c|}{Spectral problem 2} \\
\hline dof per $E$ & $E_{of}(v)$ & $E_{of}(p)$ & $E_{os}(v)$ &
$E_{os}(p)$ & $E_{of}(v)$ & $E_{of}(p)$ & $E_{os}(v)$ & $E_{os}(p)$
\rule{0pt}{9pt}
\\  \hline
1 & 0.1880 & 0.0616 & 0.1884 & 0.0405 & 0.1487 & 0.0636 & 0.1490 & 0.0428 \\
2 & 0.0427 & 0.0481 & 0.0425 & 0.0020 & 0.0833 & 0.0522 & 0.0833 & 0.0211 \\
3 & 0.0210 & 0.0481 & 0.0203 & 4.48e-04 & 0.0528 & 0.0490 & 0.0527 & 0.0099 \\
4 & 0.0107 & 0.0481 & 0.0092 & 9.35e-05 & 0.0272 & 0.0482 & 0.0267 & 0.0027 \\
5 & 0.0054 & 0.0481 & 6.57e-11 & 1.12e-14 & 0.0054 & 0.0481 & 1.45e-11 & 5.72e-15 \\
10 & 0.0054 & 0.0481 & 4.21e-12 & 7.39e-15 & 0.0054 & 0.0481 & 8.20e-12 & 7.56e-14 \\
\hline
\end{tabular}
\end{table}

%We plot the eigenvalues for both of the spectral problems in Figure
%\ref{fig:lambda}. Note that for spectral problem 2 we plot the
%eigenvalue in a descending order since we are picking the largest
%eigenvalues.

\begin{figure}[htp!]
\centering \subfigure[Spectral problem 1]{
\includegraphics[scale=.3]{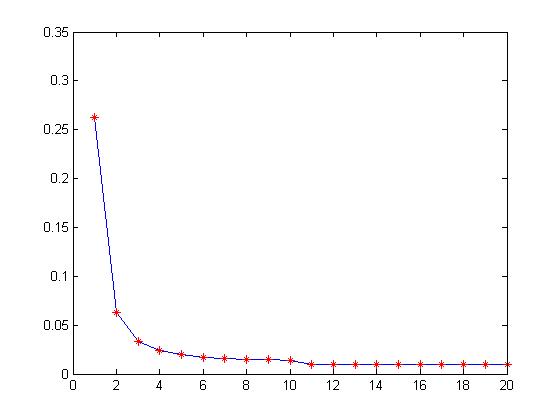}}
\subfigure[Spectral problem 2]{
\includegraphics[scale=.3]{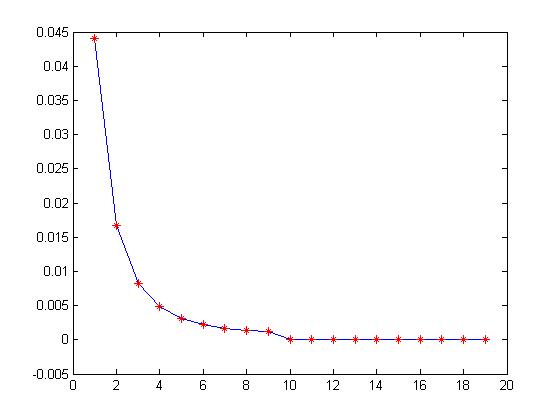}}
  \caption{Inverse of eigenvalue ($1/\lambda$) behavior for the two spectral problems.}\label{fig:lambda}
\end{figure}

In {\sc Fig}.~\ref{fig:lambda}, we show the reciprocals of the
eigenvalues for case $1$ for the permeability field $\kappa_1$ and
for a particular coarse-grid block. We also show the eigenvalue
behavior for both spectral problems. From these figures, we see that
the eigenvalues have a very sharp decay for the first $10$
eigenvalues; and this behavior corresponds to the rapid decay in the
solution errors shown in {\sc Table} \ref{TB1} and {\sc Table}
\ref{TB3}. Starting at the $11$th eigenvalue, there is no decay any
more. This situation signifies that we do not need any additional
basis function. In particular, the first $11$ eigenfunctions are
enough to achieve a machine precision spectral error, as confirmed
in {\sc Tables} \ref{TB1} and \ref{TB3}. We observe a very good
correlation (0.99) between the error and the eigenvalue behavior.

In order to see the performance of the postprocessing technique
discussed in Section \ref{sec:post}, we repeat the experiments
corresponding to {\sc Table} \ref{TB1} and compute the postprocessed
velocity, denoted as $v_p$. We define $E_{pf}(v) = \| v_p - v_f
\|_{\kappa^{-1},D} / \| v_f \|_{\kappa^{-1},D}$. The numerical
results are shown in {\sc Table} \ref{TB6}. From these results, we
clearly see that the postprocessed velocity is much more accurate
than the velocity without postprocessing.
%but with the postprocessing, see Table \ref{TB6}. In
%fact, we observe that the relative error compared with the fine
%solution $E_{of}(v)$ is much smaller in the case with
%postprocessing.
\begin{table}[h!]
\footnotesize \caption{Comparison of velocity and postprocessed
velocity, $\kappa = \kappa_1$, $n = 200$ and $N = 10$}\label{TB6}
\centering
\begin{tabular}{|c||c|c|c|c|} \hline
& \multicolumn{2}{c|}{Spectral problem 1} & \multicolumn{2}{c|}{Spectral problem 2} \\
\hline dof per $E$ & $E_{of}(v)$ & $E_{pf}(v)$  & $E_{of}(v)$ &
$E_{pf}(v)$ \rule{0pt}{9pt}
\\  \hline
1 & 0.1331 & 0.1327 & 0.1523 & 0.1525 \\
3 & 0.0569 & 0.0536 & 0.0840 & 0.0823 \\
5 & 0.0308 & 0.0232 & 0.0391 & 0.0338 \\
7 & 0.0236 & 0.0118 & 0.0278 & 0.0190 \\
9 & 0.0210 & 0.0046 & 0.0234 & 0.0114  \\
11 & 0.0208 & 0.0037 & 0.0208 & 0.0037  \\
20 & 0.0208 & 0.0037 & 0.0208 & 0.0037  \\
\hline
\end{tabular}
\end{table}

We remark that one can also consider using the curl of the velocity
in constructing the offline space. We have studied an offline space
construction that uses
\begin{equation}
a_i(v,w) = \int_{\omega_i} \text{curl}(\kappa^{-1}v)
\text{curl}(\kappa^{-1}w), \quad s_i(v,w) = \int_{\omega_i}
\kappa^{-1} v \cdot w.
\end{equation}
{\sc Table} \ref{TB7} shows the convergence of the numerical
solution obtained by using this spectral problem. As observed, the
numerical results are not as good as those shown earlier for
velocity error and for small number of basis functions.

\begin{table}[h!]
\footnotesize \caption{Convergence of the offline solution using the
curl-based spectral problem, $\kappa = \kappa_1$, $n = 200$ and $N =
10$}\label{TB7} \centering
\begin{tabular}{|c||c|c|c|c|}
\hline dof per $E$ & $E_{of}(v)$ & $E_{of}(p)$ & $E_{os}(v)$ &
$E_{os}(p)$ \rule{0pt}{9pt}
\\  \hline
1 & 0.1523 & 0.1018 & 0.1525 & 0.0519\\
3 & 0.1062 & 0.0994 & 0.1052 & 0.0447\\
5 & 0.0996 & 0.0964 & 0.0984 & 0.0373\\
7 & 0.0620 & 0.0902 & 0.0590 & 0.0108\\
9 & 0.0367 & 0.0898 &  0.0305 & 0.0024\\
11 & 0.0312 & 0.0898 & 0.0235 & 0.0013\\
20 & 0.0208 & 0.0898 & 3.90e-13 & 5.54e-15\\
\hline
\end{tabular}
\end{table}

\subsection{Oversampling technique}

Our first numerical example uses periodic coefficients. Our main
objective is to show that oversampling technique can  identify the
first-order corrector part of the solution and avoid boundary
effects. We consider the coefficient
\begin{equation*}
\kappa_{per}(x_1,x_2) = \left\{ \begin{array}{cl}
         1+\Gamma(x_1,x_2) \Pi_{i=1}^2 (0.4-|x_i-0.5|), & \mbox{if $ (x_1,x_2)\in[0.1,0.9]^2$},\\
        1, & \mbox{otherwise},\end{array} \right.
\end{equation*}
where
$$
\Gamma(x_1,x_2) = \frac{2+1.8\sin(2\pi x_1/\epsilon)}{2+1.8\sin(2\pi
x_2/\epsilon)}+\frac{2+1.8\sin(2\pi x_1/\epsilon)}{2+1.8\cos(2\pi
x_2/\epsilon)}.
$$
We consider 4 cases. {\it Case 1.} Use oversampling technique to
construct the snapshot space. When constructing the snapshot space,
we select the eigenvectors corresponding to the first $l_i^+$
eigenvalues on each coarse edge and use these eigenvectors as our
offline space. {\it Case 2.} Use oversampling technique to construct
the snapshot space. When constructing the snapshot space, we select
the eigenvectors corresponding to the first 3 eigenvalues on each
coarse edge and perform spectral problem 1 (see Section
\ref{sec:offline}) on this snapshot space
%these 3 eigenvectors in each coarse edge,
and select the eigenvectors corresponding to first $l_i$ eigenvalues
as our offline space. {\it Case 3.} Construct the snapshot space
without oversampling technique. In this case, we perform spectral
problem 1 and select the eigenvectors corresponding to first $l_i$
eigenvalues as our offline space. {\it Case 4.} Construct the
snapshot space without oversampling technique. In this case, we
perform spectral problem 2 and select the eigenvectors corresponding
to first $l_i$ eigenvalues as our offline space.
%We take $s_i$ in spectral problem 1 as follows:
%$$
%s_i(v,w) = \int_{\omega_i}\kappa^{-1}v\cdot w.
%$$
%\\
Our numerical results presented in {\sc Table} \ref{TB_over1} show
that oversampling technique does give a better performance compared
without oversampling, in general. Besides, we can obtain a much
smaller snapshot space using oversampling technique while the
accuracy of the solution is similar (see cases 2 and 3).
\begin{table}[h!]
\footnotesize \caption{Comparison of the 4 cases (relative velocity
error w.r.t. fine scale solution), $\kappa = \kappa_{per}$, $n =
200$, $N = 10$}\label{TB_over1} \centering
\begin{tabular}{|c||c|c|c|c|}
\hline dof per $E$ & Case 1 & Case 2 & Case 3 & Case 4
\\  \hline
1 & 0.0882 & 0.0985 & 0.0987 & 0.2861  \\
2 & 0.0241 & 0.0192 & 0.0206 & 0.0214  \\
3 & 0.0189 & 0.0189 & 0.0204 & 0.0210  \\
\hline
\end{tabular}
\end{table}

Next, we consider the high contrast permeability field $\kappa_1$
and compare to the previous results, see {\sc Table} \ref{TB_over2}.
Again, we see that the error is reduced if we apply oversampling
technique and the oversampling allows obtaining a small dimensional
snapshot space.
\begin{table}[h!]
\footnotesize \caption{Comparison of the 4 cases (relative velocity
error w.r.t. fine scale solution), $\kappa = \kappa_1$, $n = 200$,
$N=10$}\label{TB_over2} \centering
\begin{tabular}{|c||c|c|c|c|}
\hline dof per $E$ & Case 1 & Case 2 & Case 3 & Case 4
\\  \hline
1 & 0.1336 & 0.1332 & 0.1331 & 0.7640  \\
2 & 0.0400 & 0.0920 & 0.0916 & 0.0991  \\
3 & 0.0234 & 0.0234 & 0.0569 & 0.0593  \\
\hline
\end{tabular}
\end{table}

\subsection{Single-phase flow and transport}

We will now consider simulating single-phase flow and transport
problems by the mixed GMsFEM with spectral problem 1. Specifically,
we consider flow with zero Neumann boundary condition
\begin{equation*}
\begin{split}
-\kappa\nabla p &= {v}, \;\quad \mbox{ in } D, \\
\text{div} \, {v} &= f, \;\quad \mbox{ in } D, \\
v\cdot n &= 0, \;\quad \mbox{ on }\partial D.
\end{split}
\end{equation*}
In addition, the saturation equation is given by
\begin{equation*}
S_t + v \cdot \nabla S = r,
\end{equation*}
where $S$ is the saturation and $r$ is the source.
%Here, the velocity $v$ is solution of the above Neumann problem.
The above flow equation is solved by the mixed GMsFEM, and the
saturation equation is solved on the fine grid by the finite volume
method. Let $S^n_i$ be the value of $S$ on the fine element $\tau_i$
at time $t_n$, where $t_n = t_0 + n \Delta t$, $t_0$ is the initial
time and $\Delta t$ is the time step size chosen according to CFL
condition. Then, $S^n_i$ satisfies
\begin{equation}
|\tau_i | \frac{ S_i^{n+1} - S_i^n}{\Delta t} + \int_{\partial
\tau_i} \hat{S}^n (v\cdot n) = r_i |\tau_i |,
\label{eq:1phase_discrete}
\end{equation}
where $r_i$ is the average value of $r$ on $\tau_i$ and $\hat{S}^n$
is the upwind flux.
%That is, if
%$v\cdot n > 0 $ , then we take $\hat{S}^n = S_i^n$, otherwise take
%the value of $S$ from the other element having same edge.

In our simulations, we will take $f$ to be zero except for the
top-left and bottom-right fine-grid elements, where $f$ takes the
values of $1$ and $-1$, respectively. Moreover, we set the initial
value of $S$ to be zero. For the source $r$, we also take it as zero
except for the top-left fine element where $r=1$.
%This corresponds an injection well at the top-left corner.

%Simulation set up: we take $f$ equals $1$ on the top-left fine
%element and equals $-1$ on the bottom-right fine element. The value
%of $f$ is zero on other fine elements. The initial value of $S$ is
%zero. Also, we take $g=1$ on the top-left fine element and zero
%otherwise. This is like putting an injection well at the top-left
%corner.

\begin{figure}[htp!]
\centering \subfigure[$t = 1000$]{
\includegraphics[scale=.21]{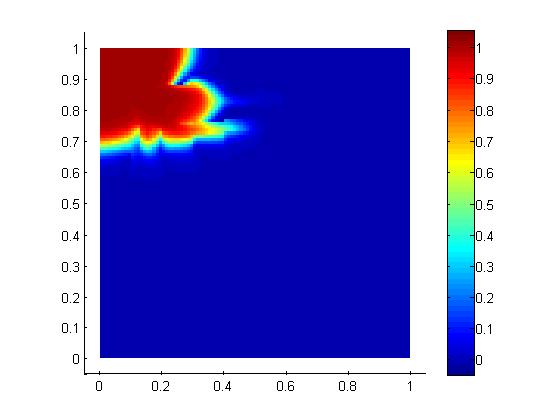}}
\subfigure[$t = 3000$]{
\includegraphics[scale=.21]{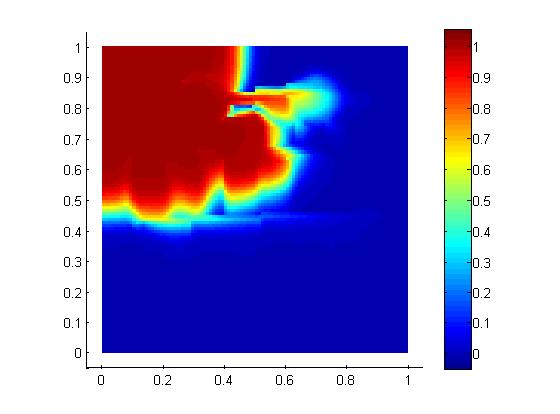}}
\subfigure[$t = 5000$]{
\includegraphics[scale=.21]{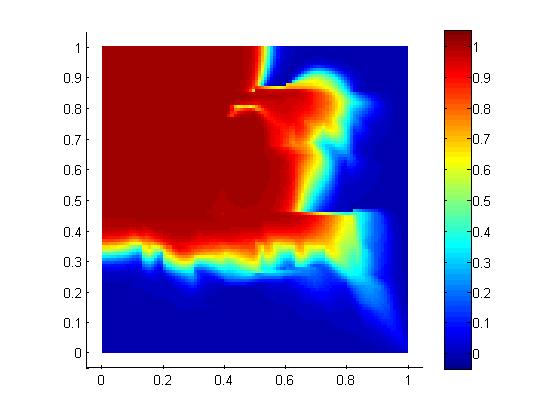}}
  \caption{Saturation solution obtained by using $v_f$ in (\ref{eq:1phase_discrete})}\label{fig:fines}
\end{figure}

\begin{figure}[htp!]
\centering \subfigure[Relative $L^2$ error = 9.0\%]{
\includegraphics[scale=.21]{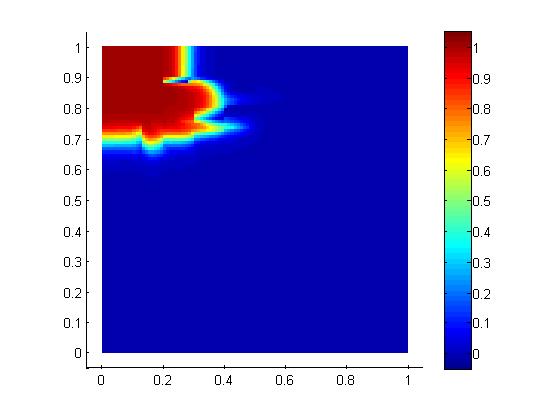}}
\subfigure[Relative $L^2$ error = 6.4\%]{
\includegraphics[scale=.21]{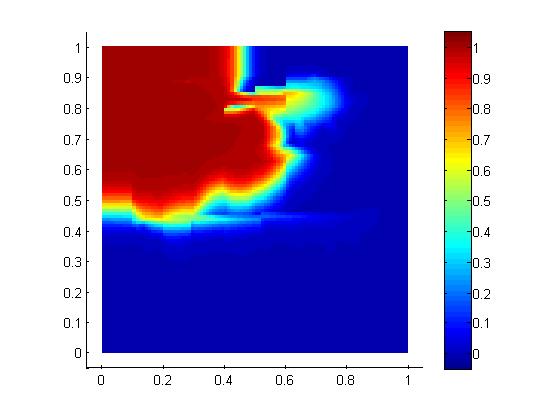}}
\subfigure[Relative $L^2$ error = 4.4\%]{
\includegraphics[scale=.21]{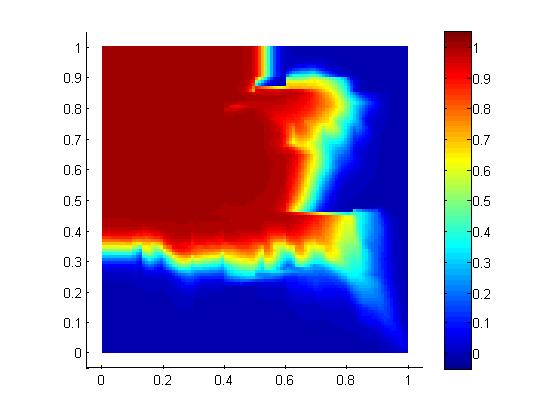}}
  \caption{Saturation solution obtained by using $v_o$ ($10\times10$ coarse grid, 1 basis per coarse edge) in (\ref{eq:1phase_discrete})}\label{fig:10off2s}
\end{figure}

\begin{figure}[htp!]
\centering \subfigure[Relative $L^2$ error = 2.0\%]{
\includegraphics[scale=.21]{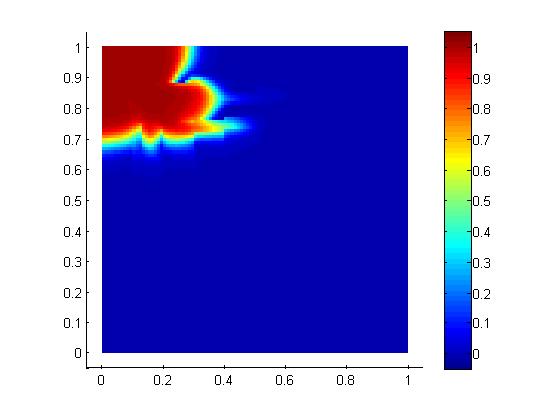}}
\subfigure[Relative $L^2$ error = 1.3\%]{
\includegraphics[scale=.21]{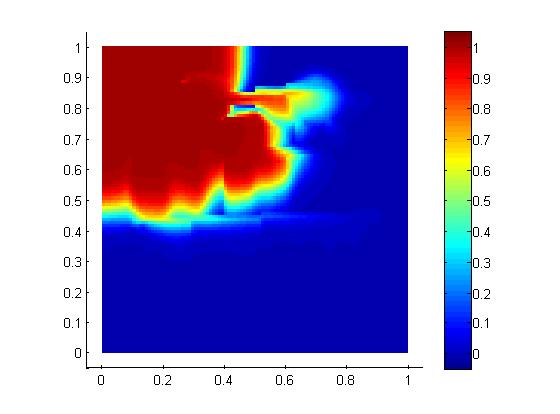}}
\subfigure[Relative $L^2$ error = 0.8\%]{
\includegraphics[scale=.21]{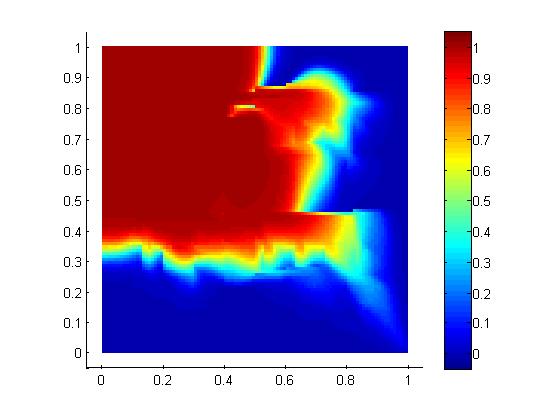}}
  \caption{Saturation solution obtained by using $v_o$ ($10\times10$ coarse grid, 3 basis per coarse edge) in (\ref{eq:1phase_discrete})}\label{fig:10off4s}
\end{figure}

\begin{figure}[htp!]
\centering \subfigure[Relative $L^2$ error = 2.0\%]{
\includegraphics[scale=.21]{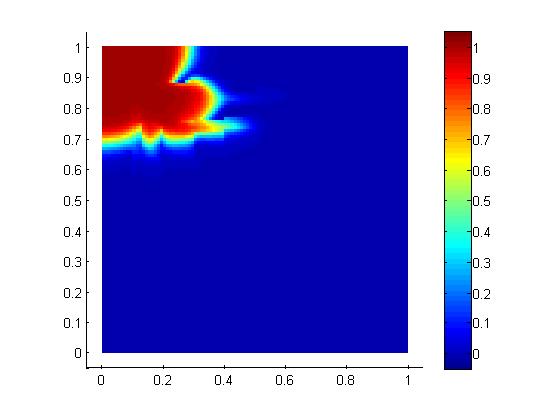}}
\subfigure[Relative $L^2$ error = 0.8\%]{
\includegraphics[scale=.21]{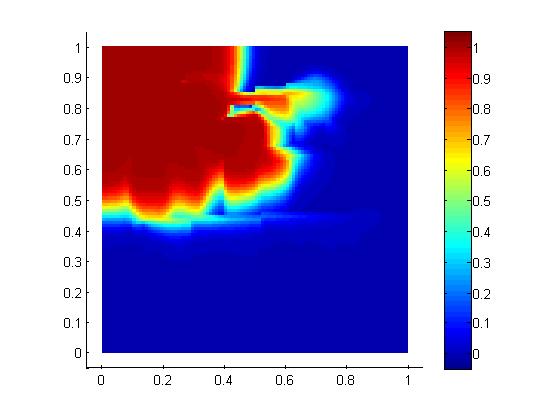}}
\subfigure[Relative $L^2$ error = 0.5\%]{
\includegraphics[scale=.21]{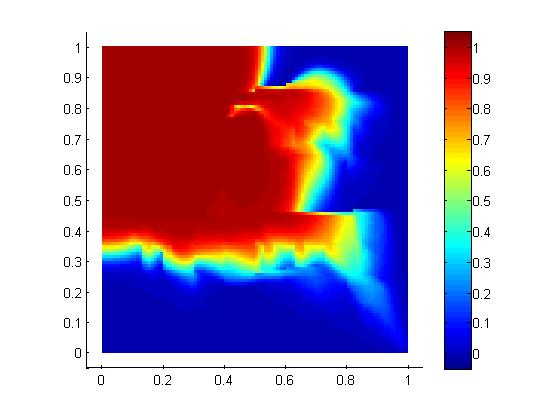}}
  \caption{Saturation solution obtained by using $v_o$ ($10\times10$ coarse grid, 5 basis per coarse edge) in (\ref{eq:1phase_discrete})}\label{fig:10off6s}
\end{figure}

In {\sc Figs}.~\ref{fig:fines}-\ref{fig:10off6s}, the saturation
plots, shown from left to right, refer to the simulations at
different times, namely; $t=1000, 3000$, and $5000$. The saturation
plots in {\sc Fig}.~\ref{fig:fines} are obtained by using the
fine-scale velocity $v_f$ in (\ref{eq:1phase_discrete}). We denote
these saturations $S_f$. Similarly, the saturation plots in {\sc
Figs}. \ref{fig:10off2s}-\ref{fig:10off6s} are obtained by using the
multiscale velocity $v_o$ in (\ref{eq:1phase_discrete}). We denote
these saturations $S_o$.
%The multiscale velocity $v_o$ is computed
%using the generalized mixed multiscale finite element method
%developed in this paper.
When selecting the multiscale basis functions, we use the first
spectral problem (\ref{eq:pod1}).
%for all the single-phase and two-phase flow simulations.
In order to see the effect of using a different number of multiscale
basis functions on each coarse edge, we repeat the simulation with
different settings. In the figures, the relative $L^2$ error refers
to the relative $L^2$ error of the saturation. We compute this as
\begin{equation*}\text{Relative error}=
\frac{\|S_o-S_f\|_{L^2(\Omega)}}{ \|S_f\|_{L^2(\Omega)}}.
\end{equation*}
In addition, we use a $10\times 10$ coarse grid for all simulations.
%Of course, when computing the relative error, we take $S_f$ and
%$S_o$ at the same simulation time.

%The coarse grid is $10\times 10$ while the fine grid is $200\times
%200$.
%With the $10\times 10$ coarse grid,
From {\sc Fig}.~\ref{fig:10off2s}, we see that if only one
multiscale basis functions are used on each coarse edge, the
relative $L^2$ error of the saturation is about $4\%$ to $9\%$. Note
that, in this case, the dimension of the velocity space
$V^0_{\text{off}}$ is only about $0.5\%$ of that of the fine scale
velocity space $V^0_h$. This shows that the mixed generalized
multiscale finite element space has a very good approximation
property. We can further reduce the relative error of saturation by
using more basis functions per coarse edge. In {\sc
Figs}.~\ref{fig:10off4s} and \ref{fig:10off6s}, we present the
relative errors for saturation when $3$ and $5$ basis functions are
used per edge respectively. We see that the errors are reduced to
approximately $2\%$. In these cases, the dimensions of the velocity
space $V^0_{\text{off}}$ are increased slightly to $1\%$ and $1.4\%$
of the fine scale velocity space $V_h^0$, respectively.



\subsection{Two-phase flow and transport}
Finally, we present our simulation results for two-phase flow and
transport problems. Consider the flow problem with zero Neumann
boundary condition
\begin{equation*}
\begin{split}
-\eta(S)\kappa\nabla p &= {v}, \;\quad \mbox{ in }D\\
\text{div} \, {v} &= f, \;\quad \mbox{ in }D \\
v\cdot n &= 0, \;\quad \mbox{ on }\partial D,
\end{split}
\end{equation*}
where
\begin{equation*}
\eta(S) = \frac{\kappa_{rw}(S)}{\mu_w}+\frac{\kappa_{ro}(S)}{\mu_o}
\end{equation*}
and
\begin{equation*}
\kappa_{rw}(S) = S^2, \;\; \kappa_{ro}(S)= (1-S)^2, \;\; \mu_w = 1,
\;\; \mu_o = 5.
\end{equation*}
The saturation equation is given by
\begin{equation*}
S_t + v \cdot \nabla F(S) = r,
\end{equation*}
where
\begin{equation*}
F(S) =
\frac{\kappa_{rw}(S)/\mu_w}{\kappa_{rw}(S)/\mu_w+\kappa_{ro}(S)/\mu_o}.
\end{equation*}
Adopting the same notations as in the single-phase flow case, we use
the following discretization for saturation
\begin{equation}
|\tau_i | \frac{ S_i^{n+1} - S_i^n}{\Delta t} + \int_{\partial
\tau_i} F(\hat{S}^n) (v\cdot n) = g_i |\tau_i
|.\label{eq:2phase_discrete}
\end{equation}
The source terms $f$ and $r$ are the same as in the single-phase
case. For the construction of the offline space, we also use the
spectral problem 1.

\begin{figure}[htp!]
\centering \subfigure[$t = 1000$]{
\includegraphics[scale=.21]{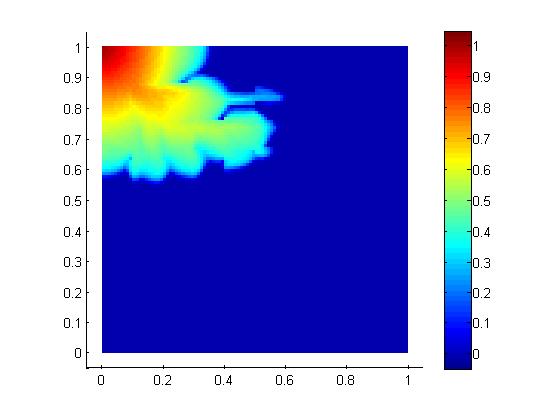}}
\subfigure[$t = 3000$]{
\includegraphics[scale=.21]{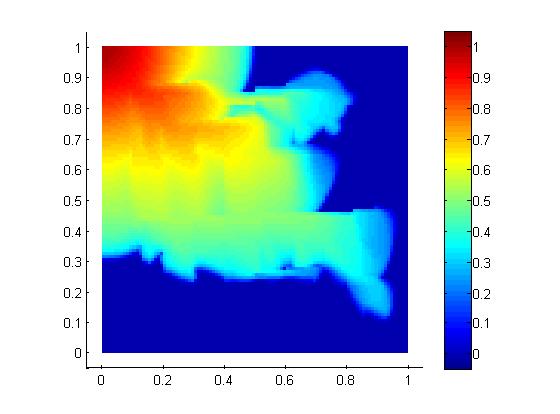}}
\subfigure[$t = 5000$]{
\includegraphics[scale=.21]{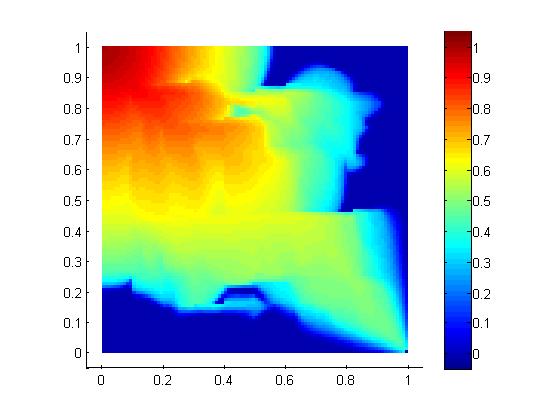}}
  \caption{Saturation solution obtained by using $v_f$ in (\ref{eq:2phase_discrete})}\label{fig:finet}
\end{figure}

\begin{figure}[htp!]
\centering \subfigure[Relative $L^2$ error = 9.3\%]{
\includegraphics[scale=.21]{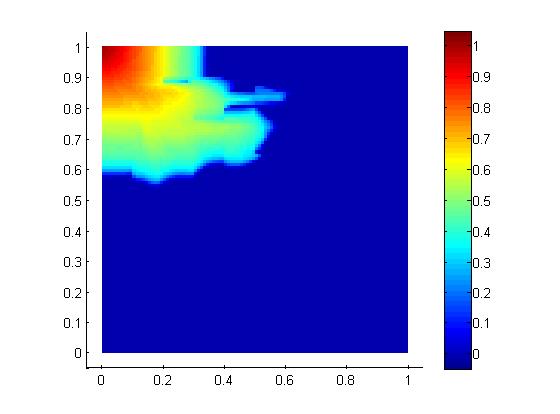}}
\subfigure[Relative $L^2$ error = 5.9\%]{
\includegraphics[scale=.21]{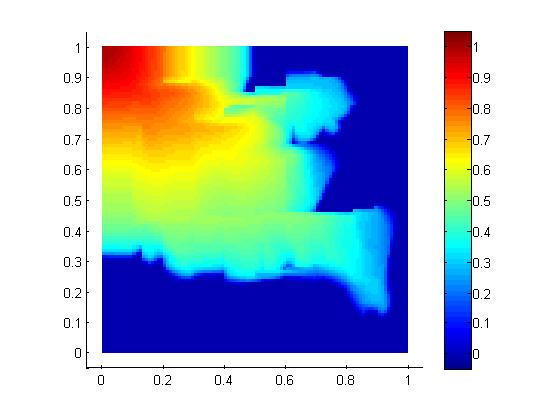}}
\subfigure[Relative $L^2$ error = 5.5\%]{
\includegraphics[scale=.21]{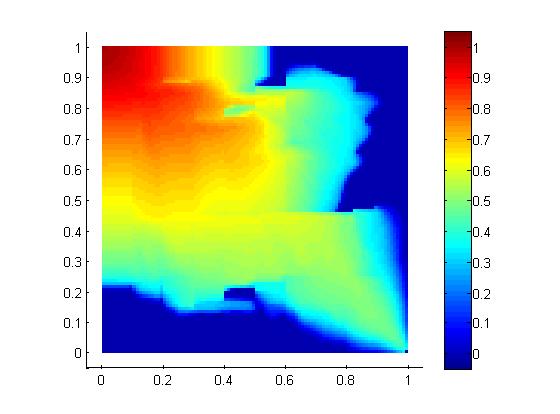}}
  \caption{Saturation solution obtained by using $v_o$ ($10\times10$ coarse grid, 1 basis per coarse edge) in (\ref{eq:2phase_discrete})}\label{fig:10off2t}
\end{figure}

\begin{figure}[htp!]
\centering \subfigure[Relative $L^2$ error = 2.8\%]{
\includegraphics[scale=.21]{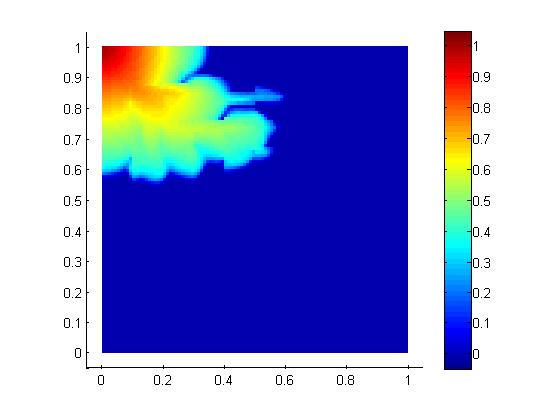}}
\subfigure[Relative $L^2$ error = 1.6\%]{
\includegraphics[scale=.21]{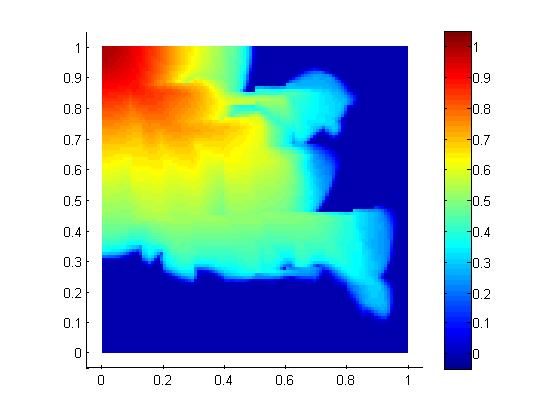}}
\subfigure[Relative $L^2$ error = 1.6\%]{
\includegraphics[scale=.21]{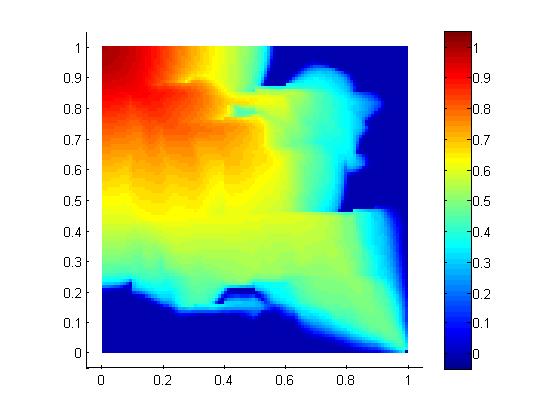}}
  \caption{Saturation solution obtained by using $v_o$ ($10\times10$ coarse grid, 3 basis per coarse edge) in (\ref{eq:2phase_discrete})}\label{fig:10off4t}
\end{figure}

\begin{figure}[htp!]
\centering \subfigure[Relative $L^2$ error = 2.6\%]{
\includegraphics[scale=.21]{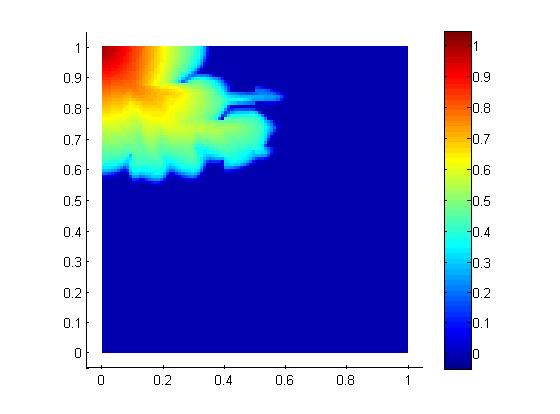}}
\subfigure[Relative $L^2$ error = 1.4\%]{
\includegraphics[scale=.21]{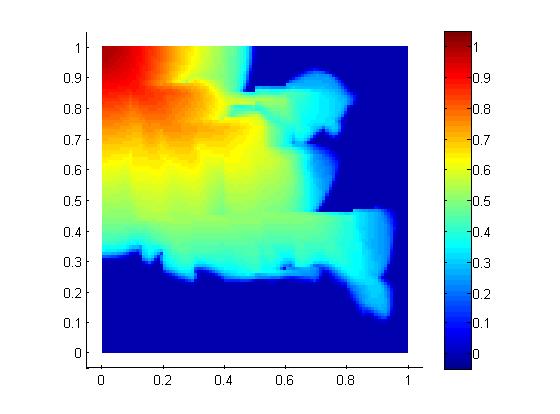}}
\subfigure[Relative $L^2$ error = 1.3\%]{
\includegraphics[scale=.21]{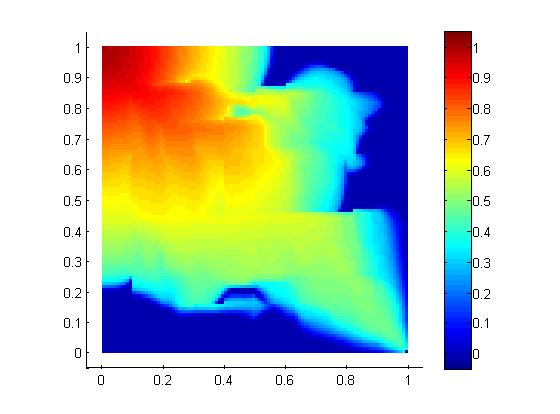}}
  \caption{Saturation solution obtained by using $v_o$ ($10\times10$ coarse grid, 5 basis per coarse edge) in (\ref{eq:2phase_discrete})}\label{fig:10off6t}
\end{figure}

In {\sc Figs}.~\ref{fig:finet}-\ref{fig:10off6t}, the saturation
plots, shown from left to right, refer to the simulations at
different times; namely, $t=1000, 3000$, and $5000$. The saturation
plots in {\sc Fig}.~\ref{fig:finet} are obtained by using the
fine-scale velocity $v_f$ in (\ref{eq:2phase_discrete}). We denote
these saturations $S_f$. Similarly, the saturation plots in {\sc
Figs}. \ref{fig:10off2t}-\ref{fig:10off6t} are obtained by using the
multiscale velocity $v_o$ in (\ref{eq:2phase_discrete}). Overall
speaking, we observe error reductions from using $1$ basis functions
per edge to $5$ basis functions per edge. In particular, for
$t=1000$, the relative error reduces from $9.3\%$ to $2.6\%$ when
using $5$ basis functions per edge, and for $t=5000$, the relative
error reduces from $5.5\%$ to $1.3\%$ when using $5$ basis functions
per edge.

%The relative $L^2$ error of the saturation is generally larger than
%the single-phase case. If we use only two basis per coarse edge, the
%relative error is about $3\%$. We also observed that in this case,
%using four basis or six basis are more or less the same as using two
%basis per coarse edge. The solution does not improve too much with
%the extra degrees of freedom beyond the first two basis in each
%coarse edge.



\begin{figure}[htp!]
\centering \subfigure[$t = 1000$]{
\includegraphics[scale=.21]{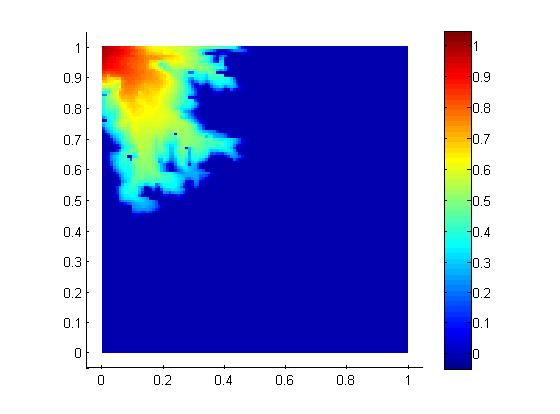}}
\subfigure[$t = 3000$]{
\includegraphics[scale=.21]{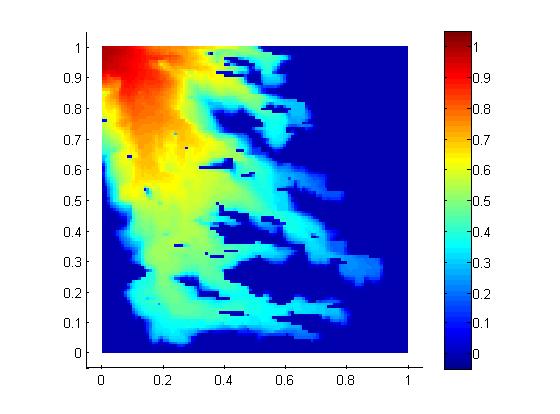}}
\subfigure[$t = 5000$]{
\includegraphics[scale=.21]{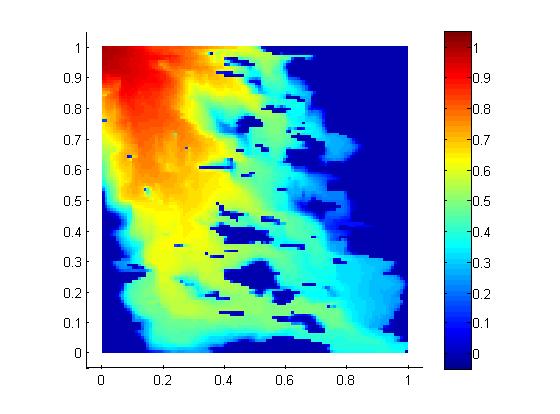}}
  \caption{Saturation solution obtained by using $v_f$ in (\ref{eq:2phase_discrete})}\label{fig:fineSPE}
\end{figure}

\begin{figure}[htp!]
\centering \subfigure[Relative $L^2$ error = 18.8\%]{
\includegraphics[scale=.21]{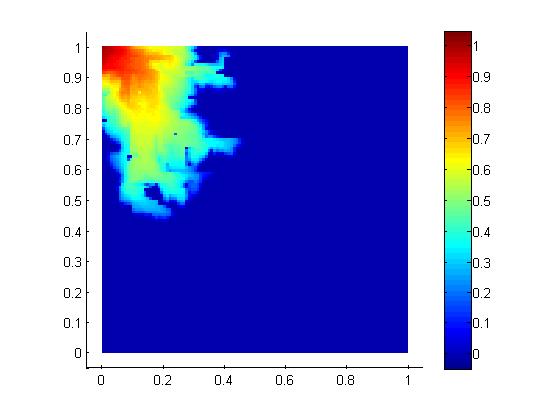}}
\subfigure[Relative $L^2$ error = 25.4\%]{
\includegraphics[scale=.21]{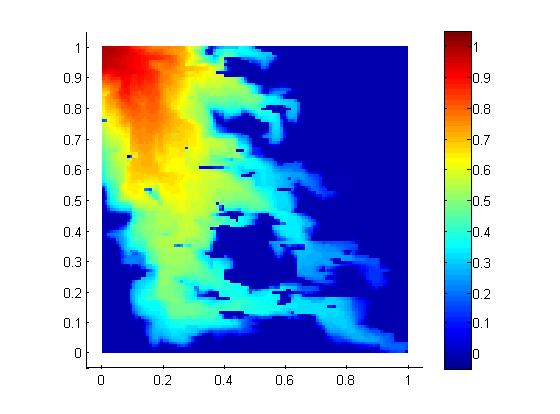}}
\subfigure[Relative $L^2$ error = 20.7\%]{
\includegraphics[scale=.21]{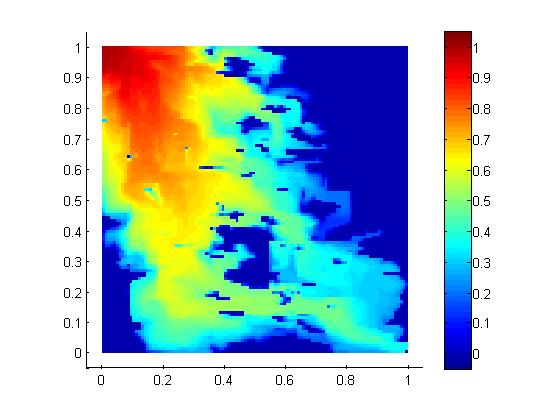}}
  \caption{Saturation solution obtained by using $v_o$ ($11\times11$ coarse grid, 1 basis per coarse edge) in (\ref{eq:2phase_discrete})}\label{fig:10off2SPE}
\end{figure}

\begin{figure}[htp!]
\centering \subfigure[Relative $L^2$ error = 5.2\%]{
\includegraphics[scale=.21]{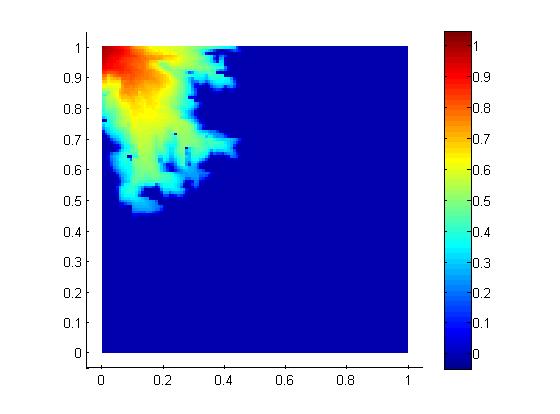}}
\subfigure[Relative $L^2$ error = 10.2\%]{
\includegraphics[scale=.21]{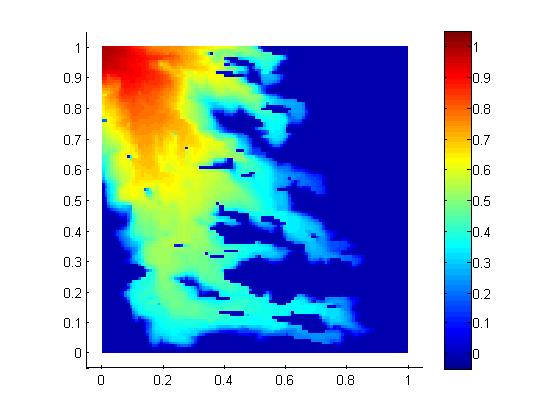}}
\subfigure[Relative $L^2$ error = 7.6\%]{
\includegraphics[scale=.21]{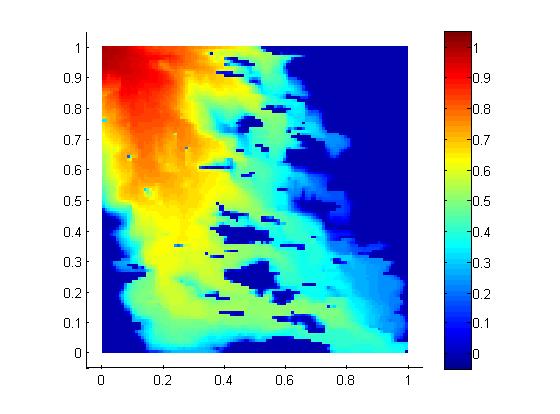}}
  \caption{Saturation solution obtained by using $v_o$ ($11\times11$ coarse grid, 3 basis per coarse edge) in (\ref{eq:2phase_discrete})}\label{fig:10off4SPE}
\end{figure}

\begin{figure}[htp!]
\centering \subfigure[Relative $L^2$ error = 3.6\%]{
\includegraphics[scale=.21]{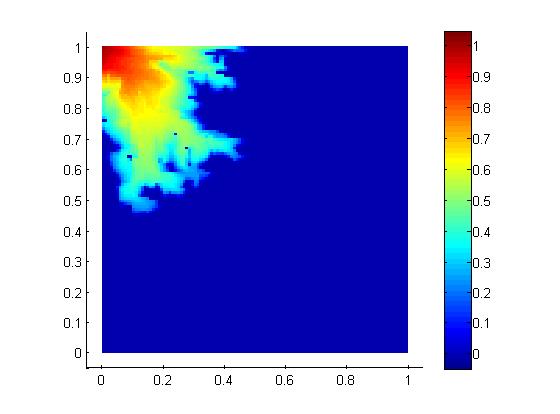}}
\subfigure[Relative $L^2$ error = 4.5\%]{
\includegraphics[scale=.21]{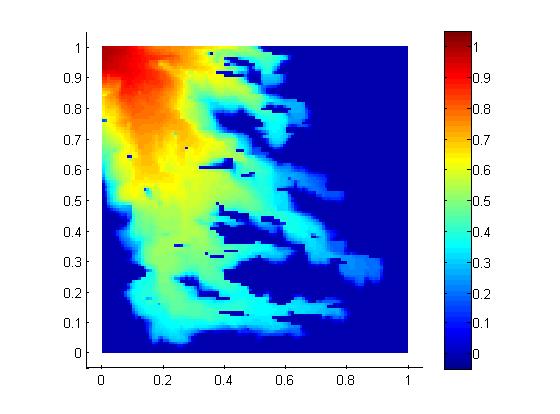}}
\subfigure[Relative $L^2$ error = 5.3\%]{
\includegraphics[scale=.21]{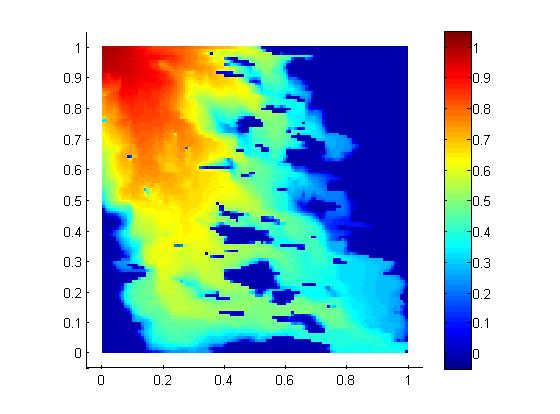}}
  \caption{Saturation solution obtained by using $v_o$ ($11\times11$ coarse grid, 5 basis per coarse edge) in (\ref{eq:2phase_discrete})}\label{fig:10off6SPE}
\end{figure}

In our last numerical example, we show the performance of our method
when applying to a more realistic permeability field. We pick the
top layer of the SPE10 permeability field (see {\sc Fig}.
\ref{fig:coef}(c)) in the following set of experiments. The model is
again the water and oil two-phase flow equations presented above.
The permeability field is originally 220 by 60, and we project it
into a fine grid of resolution 220 by 220. Then, the coarse grid is
set to be 11 by 11, which means the local grid is 10 by 10 in each
coarse block. The saturation plots are depicted in {\sc Figs}.
\ref{fig:fineSPE}-\ref{fig:10off6SPE}. In this example, we observe
that, at first glance, the multiscale saturation solution looks
similar to the fine solution if we use one multiscale basis function
per edge. However, if we take a closer look, we notice some missing
features in the water front. When we use four or six basis functions
per coarse edge, these features can be recovered correctly. This
shows the importance of these additional multiscale basis functions.
More quantitatively, we observe more error reductions from using $1$
basis functions per edge to $5$ basis functions per edge compared
with the previous examples. In particular, for $t=1000$, the
relative error reduces from $18.8\%$ to $3.6\%$ when using $5$ basis
functions per edge. Likewise, for $t=5000$, the relative error
reduces from $20.7\%$ to $5.3\%$ when using $5$ basis functions per
edge.



\section{Conclusions}

In this paper, we studied the mixed GMsFEM for constructing a mass
conservative solution of the flow equation and investigated
applications to two-phase flow and transport. The novelty of our
work is in constructing a systematic enrichment for multiscale basis
functions for the velocity field. In particular, constructing the
snapshot and the offline spaces is one of our novel contributions.
We analyze the convergence of the method and give an alternative view of eigenvalue
construction.
%Moreover, we discuss how the proposed approach can
%recover classical multiscale approaches when the problem has scale
%separation.
Besides, we study oversampling techniques and construct
snapshot vectors as the local
solutions in larger regions.
 The oversampling allows obtaining a much smaller dimensional
snapshot space and can help to improve the accuracy of the mixed GMsFEM.
Oversampling technique can be particularly helpful for problems
with scale separation.
We present numerical results and applications to single
and two-phase incompressible flow to show the performance of our method.

%\section*{References}
\bibliographystyle{plain}
\bibliography{references1}
\end{document}